\documentclass[final, 3p, times]{elsarticle}

\usepackage{enumerate}
\usepackage{graphicx}
\usepackage[dvips]{color}
\usepackage{hyperref}

\usepackage{amsmath, amsthm, amsfonts, mathrsfs, amssymb}
\usepackage{indentfirst} 
\usepackage{bm} 
\usepackage{color} 

\numberwithin{equation}{section}
\newtheorem{theorem}{Theorem}[section]

\newtheorem{lemma}{Lemma}[section]

\journal{XXX}

\begin{document}

\begin{frontmatter}

\title{On the stability and  exponential decay of the 3D MHD system with mixed partial dissipation near a equilibrium state}
\author[label1]{Xuemin Deng}
\ead{dxm@smail.xtu.edu.cn}
\author[label1]{Yuelong Xiao}
\ead{xyl@xtu.edu.cn}
\author[label2]{Aibin Zang \corref{cor1}} 
\cortext[cor1]{Corresponding author}
\ead{abzang@jxycu.edu.cn}

\address[label1]{Hunan Key Laboratory for Computation and Simulation in Science and Engineering, School of Mathematics and Computational Science, Xiangtan University, Xiangtan, Hunan, 411100, China}
\address[label2]{The Center of Applied Mathematics, Yichun University Yichun, Jiangxi, 336000, China}


\begin{abstract}

A main result of this paper establishes the global stability of the 3D MHD equations with mixed partial dissipation near a background magnetic field in the domain $\Omega=\mathbb{T}^2\times\mathbb{R}$ with $\mathbb{T}^2=[0, 1]^2$. More precisely, each velocity equation lacks its own directional dissipation, and the magnetic equation lacks vertical dissipation in the MHD system.
The key point to obtain the stability result is that we decompose the solution $(u,b)$ into the zeroth horizontal mode and the non-zeroth modes and complete the desired bound with the strong Poincar\'{e} type inequalities in the treatment of several nonlinear terms. Then we focus on the large-time behavior of the solution, where the non-zeroth modes decay exponentially in $H^2$, and the solution converges to its zeroth horizontal mode.
\end{abstract}
\end{frontmatter}
\section{Introduction}
The magnetohydrodynamic (MHD) system widely governs the dynamics of electrically conducting fluids such as plasmas, liquid metals, and salt water or electrolytes(see, e.g., \cite{BiskampDieter1993Nonlinear,Davidson2001AnIntroduction,Pippard1989inMetals}). The paper focuses on the 3D MHD system with mixed partial dissipation as follow
\begin{eqnarray}\label{SMHE}
\begin{cases}
\partial_{t} U+ U \cdot\nabla U=\nu\begin{bmatrix} \partial_{22} + \partial_{33}\\ \partial_{11} + \partial_{33}\\ \partial_{11} + \partial_{22} \end{bmatrix}U-\nabla P+B \cdot\nabla B, \hspace{0.5cm} (t, x)\in \mathbb{R}_{+}\times\Omega,\\
\partial_{t} B+ U \cdot\nabla B=\mu\left(\partial_{11} + \partial_{22} \right) B+B \cdot\nabla U, \\
div ~U=div ~B=0, \\
(U, B)\vert _{t=0}=(U_{in}, B_{in}),
\end{cases}
\end{eqnarray}
where $U=\left(U_1, U_2, U_3\right)(t, x)$ represents the three-component velocity field, $B=\left(B_1, B_2, B_3\right)(t, x)$ the magnetic field and $P=P(t, x)$ the pressure, respectively. $\nu$ and $\mu>0$ denote the kinematic viscosity and the magnetic diffusivity, respectively. Let $\Omega=\mathbb{T}^2\times\mathbb{R}$ be a spatial domain in $\mathbb{R}^3$ with $\mathbb{T}^2=[0, 1]^2$ being the 2D periodic box.



For the MHD equations  with full dissipation or partial dissipation, the global well-posedness
problem has recently attracted considerable attention and significant progress has been made
(see~\cite{CaoChongsheng2011horizontal, CaoChongsheng2011regularity, CaoChongsheng2014onlymagnetic,LinSuo2023Global,LinWuZhu2023GlobalH4,SermangeMichel1983MHD,Suo2022Global2D}).
Furthermore, in some physical experiments and numerical simulations, it was observed that electrically conducting fluids can be stabilized near a background magnetic field (see, e.g., \cite{Alemany1979Influence, Alexakis2011Two,Gallet2009Influence,Gallet2015Exact}).

We attempt to investigate the stability of the MHD system (\ref{SMHE})near a background magnetic field in the paper.
The equilibrium
$$ (U_e, B_e, P_e)=(0, e_1, 0)$$
is special steady solution of the system (\ref{SMHE}). The perturbations near a background magnetic field
$$u=U-U_e, \quad b=B-B_e, \quad p=P-P_e$$
satisfy
\begin{eqnarray}\label{smhe}
\begin{cases}
\partial_{t} u+ u \cdot\nabla u=\nu\begin{bmatrix} \partial_{22} + \partial_{33}\\ \partial_{11} + \partial_{33}\\ \partial_{11} + \partial_{22} \end{bmatrix} u-\nabla p+b \cdot\nabla b+\partial_1 b, \hspace{0.5cm} (t, x)\in \mathbb{R}_{+}\times\Omega,\\
\partial_{t} b+ u \cdot\nabla b~=\mu\left(\partial_{11} + \partial_{22} \right) b+b \cdot\nabla u+\partial_1 u, \\
div~u=div~b=0, \\
(u, b)\vert _{t=0}=(u_{in}, b_{in}) . \\
\end{cases}
\end{eqnarray}

One of our main results reads:
\begin{theorem}\label{theorem1t}
Consider the system $(\ref{smhe})$. Let the initial data $(u_{in}, b_{in})\in H^3(\Omega)$ and $div~u_{in}=div~b_{in}=0$ such that
$$\|(u_{in}, b_{in})\|_{H^3}\leq\epsilon$$
for some sufficiently small $\epsilon>0$.
Then there exists a unique global solution $(u, b)$ satisfying
\begin{align}\label{nengt}
\|(u, b)\|_{H^3}^2+\nu\int_0^t \left(\|(\partial_2, \partial_3) u_1\|_{H^3}^2
+\|(\partial_1, \partial_3) u_2\|_{H^3}^2+\|(\partial_1, \partial_2) u_3\|_{H^3}^2\right) d\tau
+\mu\int_0^t\|(\partial_1, \partial_2) b\|_{H^3}^2 d \tau\leq C\epsilon^2,
\end{align}
where the constant $C>0$ is independent of $\epsilon$ and $t$.
\end{theorem}

In recent years, the global stability or decay estimates of the anisotropic MHD equations near a background magnetic field have been studied in various domains. 

In the 2D case, when the MHD equations lack the magnetic diffusivity, the stability was investigated  in the whole space $\mathbb{R}^2$ by Hu and Lin \cite{Lin2014tability}.
Lin et al. \cite{Lin2020tability} solved the stability problems and decay estimates for the MHD system involve only vertical velocity dissipation and horizontal magnetic diffusion in $\mathbb{R}^2$.
Li, Wu and Xu \cite{XuXiaojing2020tability} verified the stability and large-time behavior of the MHD equations with only the vertical component equation in the velocity equation and the horizontal component equation in the magnetic equation involve full dissipation in $\mathbb{R}^2$.
Feng et al. \cite{Feng2023tabilityexponential} established the stability and exponential decay of the MHD equations which the velocity equation has only one direction dissipation and the magnetic equation has only damping in a spatial domain $\mathbb{T}\times\mathbb{R}$.
 Chen, Li and Chen \cite{Chen2023tability} obtained the stability and decay rate of the MHD equations with zero magnetic diffusivity in a flat strip domain $\mathbb{R}\times(0,1)$.
More recent work on the stability of  2D MHD system with various partial dissipation can be found in these references \cite{Boardman2020stability,ChenLinWu2023tability, GuoYana2022tability, Lai2021tability, Lai2022tability, LinFanghua2015mhd}.

In the 3D case, Wu and Zhu \cite{Wujiahong2021tability} proved the stability result of the MHD equations in the case when the velocity dissipation involve horizontal direction and the magnetic diffusion does not involve horizontal.
In the whole domain $\mathbb{R}^3$, Jian, Chen and Chen \cite{Jian2023tability} studied the stability and explicit decay rates of the MHD system with just damping in the velocity equation and the one direction dissipation in the magnetic equation.
The asymptotic stability and decay algebraically of the MHD equations with the velocity equation have only one-directional dissipation and the magnetic equation has horizontal dissipation was explored in $\mathbb{R}^3$ by Lin, Wu and Zhu \cite{Lin2023tabilityH3}.
Over the years there has been intensive work which the stability and large-time behavior of the MHD system in various
domain(see\cite{AbidiHammadi2017tability, DengWen2018tability, LinFanghua2014mhd, ShangHaifeng2022tability,  ZhengDahao2023tability}).


We intend to obtain decay estimates of the non-zeroth modes of the corresponding solution $(u,b)$ for the following MHD equations
\begin{eqnarray}\label{smhed}
\begin{cases}
\partial_{t} u+ u \cdot\nabla u=\nu\begin{bmatrix} \partial_{11} +\partial_{22} + \partial_{33}\\ \partial_{11} +\partial_{22} + \partial_{33}\\ \partial_{11} + \partial_{22} \end{bmatrix} u-\nabla p+b \cdot\nabla b+\partial_1 b, \hspace{0.5cm} (t, x)\in \mathbb{R}_{+}\times\Omega,\\
\partial_{t} b+ u \cdot\nabla b~=\mu\left(\partial_{11} + \partial_{22} \right) b+b \cdot\nabla u+\partial_1 u, \\
div~u=div~b=0, \\
(u, b)\vert _{t=0}=(u_{in}, b_{in}) . \\
\end{cases}
\end{eqnarray}
Hence, we decompose the solution $(u,b)$ into the zeroth horizontal mode and the non-zeroth modes in $\Omega$. More precisely, we define the zeroth horizontal mode by
\begin{align}\label{zeroth}
\overline{u}(x_3)=\int_{\mathbb{T}^2}u(x_1,x_2,x_3) dx_1dx_2, \qquad \overline{b}(x_3)=\int_{\mathbb{T}^2}b(x_1,x_2,x_3) dx_1dx_2
\end{align}
and set the non-zeroth modes as
\begin{align}\label{nonzeroth}
\widetilde{u}=u-\overline{u}, \qquad \widetilde{b}=b-\overline{b}.
\end{align}
Clearly $\overline{u}$ also represents horizontal average of $u$ and $\widetilde{u}$ the corresponding oscillation part. More details can see Section $2$. We take the average of the MHD system (\ref{smhed}) to get
\begin{eqnarray}\label{smhedh}
\begin{cases}
\partial_{t} \overline{u}+ \overline{u \cdot\nabla \widetilde{u}}+ \overline{u_3\partial_3 \overline{u}}=\nu\begin{bmatrix} \partial_{33}\overline{u}_1\\ \partial_{33}\overline{u}_2\\ 0 \end{bmatrix}-\begin{bmatrix} 0\\ 0\\ \partial_3 \overline{p} \end{bmatrix}+\overline{b \cdot\nabla \widetilde{b}}+ \overline{b_3\partial_3 \overline{b}}, \\
\partial_{t} \overline{b}+ \overline{u \cdot\nabla \widetilde{b}}+ \overline{u_3\partial_3 \overline{b}}=\overline{b \cdot\nabla \widetilde{u}}+ \overline{b_3\partial_3 \overline{u}}, \\
\partial_3 \overline{u}_3=\partial_3 \overline{b}_3=0, \\
(\overline{u}, \overline{b})\vert _{t=0}=(\overline{u}_{in}, \overline{b}_{in}) . \\
\end{cases}
\end{eqnarray}
Taking the difference of the system (\ref{smhed}) and the system (\ref{smhedh}), we find
\begin{eqnarray}\label{smhedv}
\begin{cases}
\partial_{t} \widetilde{u}+ \widetilde{u \cdot\nabla \widetilde{u}}+ \widetilde{u}_3\partial_3 \overline{u}=\nu\begin{bmatrix} \partial_{11} +\partial_{22} + \partial_{33}\\ \partial_{11} +\partial_{22} + \partial_{33}\\ \partial_{11} + \partial_{22} \end{bmatrix} \widetilde{u}-\nabla \overline{p} +\widetilde{b \cdot\nabla \widetilde{b}}+ \widetilde{b}_3\partial_3 \overline{b}+\partial_3 \widetilde{b}, \\
\partial_{t} \widetilde{b}+ \widetilde{u \cdot\nabla \widetilde{b}}+ \widetilde{u}_3\partial_3 \overline{b}=\mu\left(\partial_{11} + \partial_{22} \right) \widetilde{b} +\widetilde{b \cdot\nabla \widetilde{u}}+ \widetilde{b}_3\partial_3 \overline{u}+\partial_3 \widetilde{u}, \\
div~\widetilde{u}=div~\widetilde{b}=0, \\
(\widetilde{u}, \widetilde{b})\vert _{t=0}=(\widetilde{u}_{in}, \widetilde{b}_{in}) . \\
\end{cases}
\end{eqnarray}
Obviously, the system (\ref{smhed}) conforms to Theorem \ref{theorem1t}.
Now we are ready to present another result of this paper on the large-time behavior of the solution.
\begin{theorem}\label{theorem2t}
Assume that $(u, b)$ is the solution to the system (\ref{smhed}) with the initial data $(u_{in}, b_{in})$ under the conditions in Theorem \ref{theorem1t}.
Then the non-zeroth modes $(\widetilde{u}, \widetilde{b})$ of $(u, b)$ decays exponentially and satisfies
\begin{align*}
\|(\widetilde{u}, \widetilde{b})\|_{H^2}\leq \|(u_{in}, b_{in})\|_{H^2}e^{-\frac{\lambda}{2} t},
\end{align*}
where $\lambda=min\{\nu, \mu\}$. Therefore, the system (\ref{smhed}) converges to the system as follow
\begin{eqnarray}\label{smhedhl}
\begin{cases}
\partial_{t} \overline{u}_1+ \overline{u}_3\partial_3 \overline{u}_1=\nu\partial_{33}\overline{u}_1+ \overline{b}_3\partial_3 \overline{b}_1, \\
\partial_{t} \overline{u}_2+ \overline{u}_3\partial_3 \overline{u}_2=\nu\partial_{33}\overline{u}_2+ \overline{b}_3\partial_3 \overline{b}_2, \\
\partial_{t} \overline{b}_1+ \overline{u}_3\partial_3 \overline{b}_1= \overline{b}_3\partial_3 \overline{u}_1, \\
\partial_{t} \overline{b}_2+ \overline{u}_3\partial_3 \overline{b}_2= \overline{b}_3\partial_3 \overline{u}_2, \\
\overline{u}_3(x_3,t)=\overline{u}_{in,3},  \quad \overline{b}_3(x_3,t)=\overline{b}_{in,3}.\\
\end{cases}
\end{eqnarray}
\end{theorem}

The detailed statements of Theorem \ref{theorem1t} and \ref{theorem2t} can be found in Sections~\ref{sec3t} and \ref{sec4t}.
Now we will briefly explain some of the main difficulties and proof lines.
The proof of this paper is divided into twofold.
First step, we structure the energy functional
\begin{align}\label{energyfunctionalt}
E(t)=&\underset{0\leq\tau\leq t}{sup}\|(u, b)\|_{H^3}^2+\nu\int_0^t \left(\|(\partial_2, \partial_3) u_1\|_{H^3}^2
+\|(\partial_1, \partial_3) u_2\|_{H^3}^2+\|(\partial_1, \partial_2) u_3\|_{H^3}^2\right) d\tau
+\mu\int_0^t\|(\partial_1, \partial_2) b\|_{H^3}^2 d \tau
\end{align}
to satisfy
\begin{align}\label{ett}
E(t)\leq E(0)+CE(t)^{\frac{3}{2}}, \quad \forall t\geq 0.
\end{align}
Then the stability (\ref{nengt}) is obtained by using the bootstrapping argument.
However, the process is nontrivial.
For the MHD system (\ref{smhe}), we originally tried to get stability results in the whole space $\mathbb{R}^3$.
Due to the velocity equation lacks its own directional dissipation and the magnetic equation lacks Vertical diffusion, we use the following anisotropic inequality
\begin{align}\label{geR3}
\int_{\mathbf{R}^3} |fgh|dx\leq C \|f\|^{\frac{1}{2}} \|\partial_1 f\|^{\frac{1}{2}} \|g\|^{\frac{1}{2}} \|\partial_2 g\|^{\frac{1}{2}} \|h\|^{\frac{1}{2}} \|\partial_3 h\|^{\frac{1}{2}},
\end{align}
to get the hardest one
\begin{align}\label{Hardterm}
\sum_{i,k=1}^2 \int \partial_i u_i\partial_3^3 b_k\partial_3^3 b_k dx
\leq C\sum_{i,k=1}^2\|\partial_3^3 b_k\|^\frac{1}{2} \|\partial_1\partial_3^3 b_k\|^\frac{1}{2} \|\partial_3^3 b_k\|^\frac{1}{2} \|\partial_2\partial_3^3 b_k\|^\frac{1}{2} \|\partial_i u_i\|^\frac{1}{2} \|\partial_{3i} u_i\|^\frac{1}{2}.
\end{align}
The rest of the process, the divergence free condition $div~u=div~b=0$ is also invalid here. The subterm $\|\partial_3^3 b_k\|$, $\|\partial_i u_i\|$ cannot be directly controlled by
\begin{align*}
\sum_{i=1}^2 \|\partial_i b\|_{H^3},~~\sum_{i=1}^3\left(\|\partial_{i'} u_i\|_{H^3} +\|\partial_{i''} u_i\|_{H^3}\right),
\end{align*}
respectively.
Where $i, i', i''\in \{1, 2, 3\}$ satisfies $\varepsilon_{ii'i''}\neq0$ ~and~ $\varepsilon_{ii'i''}$ denotes the Ricci symbol.
We expect to control the right side of (\ref{Hardterm}) using a Poincar\'{e} inequality.
But the right side of the standard Poincar\'{e} type inequality  involves the full gradient not partial direction.
The idea is inspired by the approach used in the paper \cite{JiRuihong2023tability}.
In the domain $\Omega=\mathbb{T}^2\times\mathbb{R}$,
the non-zeroth modes $\widetilde{u}, ~\widetilde{b}$ obeys a strong version of the Poincar\'{e} type inequality
\begin{align}\label{Poincareub}
\|\widetilde{u}\|_{H^s(\Omega)}\leq C\|\partial_i\widetilde{u}\|_{H^s(\Omega)}, \quad \|\widetilde{b}\|_{H^s(\Omega)}\leq C\|\partial_i\widetilde{b}\|_{H^s(\Omega)}.
\end{align}
Using the definition of (\ref{zeroth}) and (\ref{nonzeroth}), we decompose the left side of the (\ref{Hardterm}) to get
\begin{align}\label{Hardtermt}
\sum_{i,k=1}^2 \int \partial_i u_i\partial_3^3 b_k\partial_3^3 b_k dx= \sum_{i,k=1}^2 \int \partial_i (\overline{u}_i+\widetilde{u}_i)\partial_3^3 (\overline{b}_k+\widetilde{b}_k)\partial_3^3 (\overline{b}_k+\widetilde{b}_k) dx.
\end{align}
Thanks to the anisotropic inequality (\ref{fge2}) and (\ref{fge3}) of Lemma \ref{lemget} and Lemma \ref{lempoincare} help complete the bounds.
Then we complete the proof in the inequality  (\ref{energyfunctionalt}) and then we get the uniqueness part.
The second step is to explore the large-time behavior problem in $H^2$. Due to the partial dissipation, the Fourier splitting method is ineffective.
Hence, we first obtain the decays exponentially of the non-zeroth modes solution $(\widetilde{u},\widetilde{b})$ for the MHD system (\ref{smhedv}) by using standard energy method.
Then the system (\ref{smhed}) converges to the system (\ref{smhedhl}).
We frequently use Lemma \ref{lemget} and Lemma \ref{lempoincare} in nonlinear term bound.
In particular, we introduce the following anisotropic inequality
\begin{align*}
\int |\widetilde{f}\widetilde{g}h|dx\leq& C \|\widetilde{f}\|^{\frac{1}{4}} \|\partial_i \widetilde{f}\|^\frac{1}{4}\|\partial_3 \widetilde{f}\|^\frac{1}{4}\|\partial_{3i} \widetilde{f}\|^\frac{1}{4}\|\widetilde{g}\|^{\frac{1}{2}} \|\partial_j \widetilde{g}\|^\frac{1}{2}\|h\|
\end{align*}
to ensure that the exponential decay results are in $H^2$.


The rest of the paper is organized as follows. In Section~\ref{sec2t}, some notations and definitions are introduced. We mainly present the definition of the zeroth mode $\overline{f}$ and the non-zeroth modes $\widetilde{f}$ and some related properties and anisotropic upper bounds. In Section ~\ref{sec3t}, we show the proof of theorem \ref{theorem1t}. Finally, the large-time behavior problem is analyzed in Section ~\ref{sec4t}.

\section{Preliminaries}\label{sec2t}
In this section, we will describe in detail the relevant definitions and properties of the zeroth horizontal mode and the non-zeroth modes. The various anisotropic inequalities involved in the proof are presented. Throughout this paper for simplicity, C stands for some real positive constant. Denote $$\|f\|^2:=\|f\|^2_{L^2(\Omega)}, \quad\|f\|_{L_t^2 H^s}:=\int_0^t\|f\|^2_{H^s(\Omega)} d\tau, ~s\geq0 $$
and
$$\|(\partial^{s_1}, \partial^{s_2}) (f,g)\|^2 := \|\partial^{s_1} (f,g)\|^2+\|\partial^{s_2} (f,g)\|^2:= \|\partial^{s_1} f\|^2+ \|\partial^{s_1} g\|^2+\|\partial^{s_2} f\|^2+\|\partial^{s_2} g\|^2,  ~s_1,s_2\geq0,$$
the norm equivalence
$$\|f\|_{H^3 }^2\sim \|f\|^2+\sum_{i=1}^3\|\partial_i^3f\|^2,$$
where $f=f(t)$ and $ g=g(t)$ are two measurable functions.

In the domain $\Omega=\mathbb{T}^2\times\mathbb{R}$, We recall the zeroth mode and the non-zeroth modes
\begin{align}\label{fdecom}
\overline{f}(x_3)=\int_{\mathbb{T}^2}f(x_1,x_2,x_3) dx_1x_2, \qquad \widetilde{f}=f-\overline{f}.
\end{align}

We introduce the following detailed properties of the zeroth mode and the non-zeroth modes.
\begin{lemma}(see \cite{JiRuihong2023tability}.)\label{lempro}
Let $s\geq 0$ be an integer. Suppose that $f$ is sufficiently regular and a divergence-free vector field  in $\Omega=\mathbb{T}^2\times\mathbb{R}$. Let $\overline{f}(x_3)$~and~$\widetilde{f}(x_1,x_2,x_3)$ be defined as in (\ref{fdecom}).
\begin{enumerate}[(i)]
\rm\item $\overline{\partial_i^s f}=\partial_i^s\overline{f}, \quad\widetilde{\partial_i^s f}=\partial_i^s\widetilde{f}. \quad(i=1,2,3)$
\rm\item $(\overline{f},\widetilde{f})_{H^s}=0, \quad\|f\|_{H^s}=\|\overline{f}\|_{H^s}+\|\widetilde{f}\|_{H^s}.$
\rm\item $div~\overline{f}=0, \quad div~\widetilde{f}=0. $
\end{enumerate}
\end{lemma}

The following lemma is crucial for deriving the triple product integral.
\begin{lemma}(see \cite{JiRuihong2023tability}.)\label{lem1D}
Let $f\in H^1(\mathbb{T})$ and Let $\widetilde{f}(x_1,x_2,x_3)$ be defined as in (\ref{fdecom}). Then
\begin{align}\label{f1D}
\|f\|_{L^{\infty}(\mathbb{T})}\leq &\sqrt{2}\|f\|_{L^2(\mathbb{T})}^\frac{1}{2} \left(\|f\|_{L^2(\mathbb{T})}+\|f'\|_{L^2(\mathbb{T})} \right)^\frac{1}{2} \\
\|\widetilde{f}\|_{L^{\infty}(\mathbb{T})}\leq &\sqrt{2}\|\widetilde{f}\|_{L^2(\mathbb{T})}^\frac{1}{2} \|\widetilde{f}'\|_{L^2(\mathbb{T})}^\frac{1}{2}
\end{align}
\end{lemma}
The following anisotropic inequalities are obtained by invoking the one-dimensional inequalities in lemma \ref{lem1D} and Minkowski's inequality. See reference \cite{CaoChongsheng2011regularity} for details of the proof.
\begin{lemma}\label{lemget}
Let $i=1, ~j=2$ or $i=2, ~j=1$. Suppose that $f, g, h, \partial_i f, \partial_j g $ and $\partial_3 h$ are all in $L^2(\Omega)$ with $\Omega=\mathbb{T}^2\times\mathbb{R}$. Then
\begin{align}\label{fge1}
\int |fgh|dx\leq& C \|f\|^{\frac{1}{2}} \left(\|f\|+\|\partial_1 f\| \right)^\frac{1}{2}\|g\|^{\frac{1}{2}} \left(\|g\|+\|\partial_2 g\| \right)^\frac{1}{2} \|h\|^{\frac{1}{2}} \|\partial_3 h\|^{\frac{1}{2}},
\end{align}
\begin{align}\label{fge2}
\int |\widetilde{f}g\widetilde{h}|dx\leq& C \|\widetilde{f}\|^{\frac{1}{2}} \|\partial_i \widetilde{f}\|^\frac{1}{2}\|g\|^{\frac{1}{2}} \left(\|g\|+\|\partial_j g\| \right)^\frac{1}{2}\|h\|^{\frac{1}{2}} \|\partial_3 h\|^{\frac{1}{2}},
\end{align}
and
\begin{align}\label{fge3}
\int |\widetilde{f}\widetilde{g}\widetilde{h}|dx\leq& C \|\widetilde{f}\|^{\frac{1}{2}} \|\partial_1 \widetilde{f}\|^\frac{1}{2}\|\widetilde{g}\|^{\frac{1}{2}} \|\partial_2 \widetilde{g}\|^\frac{1}{2}\|\widetilde{h}\|^{\frac{1}{2}} \|\partial_3 \widetilde{h}\|^{\frac{1}{2}}.
\end{align}
\end{lemma}
In the whole space $\mathbb{R}^3$, the version of the anisotropic inequality
\begin{align*}
\int_{\mathbb{R}^3} |fgh|dx\leq& C \|f\|^{\frac{1}{4}} \|\partial_1 f\|^\frac{1}{4}\|\partial_2 f\|^{\frac{1}{4}} \|\partial_{12} f\|^\frac{1}{4}\|g\|^{\frac{1}{2}} \|\partial_3 g\|^\frac{1}{2} \|h\|
\end{align*}
is provided by Wu et al\cite{Wujiahong2021tability}.  Similarly, we can deduce the following Lemma in the spatial domain $\mathbb{T}^2\times\mathbb{R}$.
\begin{lemma}\label{lemge3t}
Let $i=1, ~j=2$ or $i=2, ~j=1$. Suppose that $f, g, h, \partial_i f, \partial_3 f, \partial_{3i} f $ and $\partial_j g$ are all in $L^2(\Omega)$ with $\Omega=\mathbb{T}^2\times\mathbb{R}$. Then
\begin{align}\label{fge13}
\int |fgh|dx\leq& C \|f\|^{\frac{1}{4}} \left(\|f\|+\|\partial_i f\| \right)^\frac{1}{4}\|\partial_3 f\|^{\frac{1}{4}} \left(\|\partial_3 f\|+\|\partial_{3i} f\| \right)^\frac{1}{4}\|g\|^{\frac{1}{2}} \left(\|g\|+\|\partial_j g\| \right)^\frac{1}{2} \|h\|
\end{align}
and
\begin{align}\label{fge23}
\int |\widetilde{f}\widetilde{g}h|dx\leq& C \|\widetilde{f}\|^{\frac{1}{4}} \|\partial_i \widetilde{f}\|^\frac{1}{4}\|\partial_3 \widetilde{f}\|^\frac{1}{4}\|\partial_{3i} \widetilde{f}\|^\frac{1}{4}\|\widetilde{g}\|^{\frac{1}{2}} \|\partial_j \widetilde{g}\|^\frac{1}{2}\|h\|.
\end{align}
\end{lemma}

\begin{lemma}\label{lempoincare}
Let $\Omega=\mathbb{T}^2\times\mathbb{R}$ and $\partial_i\widetilde{f}\in H^s(\Omega), ~s\geq 0, ~i=1,2$. Assume that the oscillation part $\widetilde{f}(x_1,x_2,x_3)$ be defined as in (\ref{fdecom}). Then there exists $C\geq 0$ such that
\begin{align}\label{fpoincare}
\|\widetilde{f}\|_{H^s(\Omega)}\leq C\|\partial_i\widetilde{f}\|_{H^s(\Omega)}.
\end{align}
\end{lemma}

\section{Proof of Theorem \ref{theorem1t}.} \label{sec3t}
This section discloses the stability of the system (\ref{smhe}), as shown in Theorem \ref{theorem1t}.
The important technique here is bootstrap argument, which first establishes a priori estimate and then applies this method to obtain inequality (\ref{nengt}).

\begin{proof}
First, we make $L^2$ estimate.
Taking the $L^2$ inner product of the equations (\ref{smhe}) with $(u, b)$, we obtain
\begin{align*} 
\frac{1}{2}\frac{d}{dt}\|(u, b)\|^2+\nu\left(\|(\partial_2, \partial_3) u_1\|^2
+\|(\partial_1, \partial_3) u_2\|^2+\|(\partial_1, \partial_2) u_3\|^2\right)+\mu\|(\partial_1, \partial_2) b\|^2=0.
\end{align*}
Integrating in time one has 
\begin{align*}
\|(u, b)\|^2+\nu\int_0^t \|(\partial_2, \partial_3) u_1\|^2
+\|(\partial_1, \partial_3) u_2\|^2+\|(\partial_1, \partial_2) u_3\|^2 d\tau
+\mu\int_0^t\|(\partial_1, \partial_2) b\|^2 d \tau=0.
\end{align*}
Applying $\partial_i^3$ to (\ref{smhe}) and taking $L^2$ inner product with $(\partial_i^3 u, \partial_i^3 b)$ and then taking summation for $i$ from $1$ to $3$, we obtain
\begin{align*}
\frac{1}{2}\frac{d}{d t}&\|(u, b)\|_{H^3}^2+\nu\left(\|(\partial_2, \partial_3) u_1\|_{H^3}^2+\|(\partial_1, \partial_3) u_2\|_{H^3}^2+\|(\partial_1, \partial_2) u_3\|_{H^3}^2\right)+\mu\|(\partial_1, \partial_2) b\|_{H^3}^2\\
&\leq-\sum_{m,i=1}^3 C_3^m\int \partial_i^m u\cdot\nabla \partial_i^{3-m} u\cdot\partial_i^3 u dx+\sum_{m,i=1}^3  C_3^m\int \partial_i^m b\cdot\nabla \partial_i^{3-m} b\cdot\partial_i^3 u dx\\
&~~~-\sum_{m,i=1}^3 C_3^m\int \partial_i^m u\cdot\nabla \partial_i^{3-m} b\cdot\partial_i^3 b dx+\sum_{m,i=1}^3  C_3^m\int \partial_i^m b\cdot\nabla  \partial_i^{3-m}u\cdot\partial_i^3 b dx\\
&:=I+J+K+L.
\end{align*}
We notice a fact that
\begin{align*}
\sum_{i=1}^3 \int \partial_i^3 \partial_1 b\cdot \partial_i^3 u dx+\sum_{i=1}^3 \int \partial_i^3 \partial_1 u\cdot\partial_i^3 b dx =0
\end{align*}
and
\begin{align*}
\sum_{i=1}^3  \left(\int  b\cdot\nabla\partial_i^3 b\cdot\partial_i^3 u dx+\int b\cdot\nabla\partial_i^3 u\cdot\partial_i^3 b dx\right)=0.
\end{align*}
\noindent\textbf {Estimate of $I$.} We decompose it into
\begin{align*}
I=&-\sum_{m,i=1}^3 C_3^m\int \partial_i^m u_i \partial_i\partial_i^{3-m} u_i\partial_i^3 u_i dx -\sum_{m,i=1}^3\sum_{j\neq i} C_3^m\int \partial_i^m u_j \partial_j\partial_i^{3-m} u_i\partial_i^3 u_i dx\\
&-\sum_{m,i=1}^3\sum_{k\neq i} C_3^m\int \partial_i^m u_i \partial_i\partial_i^{3-m} u_k\partial_i^3 u_k dx-\sum_{m,i=1}^3\sum_{j\neq i,k\neq i} C_3^m\int \partial_i^m u_j \partial_j\partial_i^{3-m} u_k\partial_i^3 u_k dx\\
:=&I_1+I_2+I_3+I_4.
\end{align*}
For the first item $I_1$, we decompose it into three terms
\begin{align*}
I_1=&-\sum_{m=1}^3 C_3^m\int \partial_1^m u_1 \partial_1\partial_1^{3-m} u_1\partial_1^3 u_1 dx-\sum_{m=1}^3 C_3^m\int \partial_2^m u_2 \partial_2\partial_2^{3-m} u_2\partial_2^3 u_2 dx\\&-\sum_{m=1}^3 C_3^m\int \partial_3^m u_3 \partial_3\partial_3^{3-m} u_3\partial_3^3 u_3 dx:=I_{11}+I_{12}+I_{13}.
\end{align*}
To estimate the term $I_{11}$, we rewrite
\begin{align*}
I_{11}=&-3\sum_{j,k=2}^3 \int \partial_1 u_1 \partial_1^2\partial_j u_j\partial_1^2\partial_k u_k dx-\sum_{m,j,k=2}^3 C_3^m\int \partial_1^{m-1}\partial_j u_j \partial_1\partial_1^{3-m} u_1\partial_1^2\partial_k u_k dx:=I_{111}+I_{112}.
\end{align*}
Due to the divergence free condition $div~u=0$, we can get the following bounded by the inequality (\ref{fge1}) and Young inequality.
\begin{align*}
I_{111}\leq&C\sum_{j,k=2}^3\|\partial_1^2\partial_j u_j\|^\frac{1}{2} \left(\|\partial_1^2\partial_j u_j\|+\|\partial_1^3\partial_j u_j\|\right)^\frac{1}{2} \|\partial_1^2\partial_k u_k\|^\frac{1}{2} \left(\|\partial_1^2\partial_k u_k\|+\|\partial_1^3\partial_k u_k\|\right)^\frac{1}{2} \|\partial_1 u_1\|^\frac{1}{2} \|\partial_{31} u_1\|^\frac{1}{2}\\
\leq&C\|u\|_{H^3} \left(\|\partial_1 u_2\|_{H^3}^2+\|\partial_1 u_3\|_{H^3}^2\right)
\end{align*}
and
\begin{align*}
I_{112}\leq&C\sum_{m,j,k=2}^3\|\partial_1^2\partial_k u_k\|^\frac{1}{2} \left(\|\partial_1^2\partial_k u_k\|+\|\partial_1^3\partial_k u_k\|\right)^\frac{1}{2} \|\partial_1^{m-1}\partial_j u_j\|^\frac{1}{2} \left(\|\partial_1^{m-1}\partial_j u_j\|+\|\partial_{2j}\partial_1^{m-1} u_j\|\right)^\frac{1}{2} \|\partial_1^{4-m} u_1\|^\frac{1}{2} \|\partial_3\partial_1^{4-m} u_1\|^\frac{1}{2}\\
\leq&C\|u\|_{H^3} \left(\|\partial_1 u_2\|_{H^3}^2+\|\partial_1 u_3\|_{H^3}^2\right).
\end{align*}
Similarly,
\begin{align*}
I_{12}\leq C\|u\|_{H^3} \left( \|\partial_2 u_1\|_{H^3}^2+ \|\partial_2 u_3\|_{H^3}^2\right), \hspace{0.2cm}
I_{13}\leq C\|u\|_{H^3} \left( \|\partial_3 u_1\|_{H^3}^2 +\|\partial_3 u_2\|_{H^3}^2\right).\\
\end{align*}
 For the term $I_2$, $j$ can choose $i'$ or $i''$. By using the divergence free condition $div~u=0$, the term $I_2$ can be expanded under the form
\begin{align*}
I_2=&-\sum_{i,m=1}^3\sum_{j\neq i} C_3^m\int \partial_i^m u_j \partial_j\partial_i^{3-m} u_i \partial_i^2\partial_{i'} u_{i'} dx-\sum_{i,m=1}^3\sum_{j\neq i} C_3^m\int \partial_i^m u_j \partial_j\partial_i^{3-m} u_i \partial_i^2\partial_{i''} u_{i''} dx:=I_{21}+I_{22}.
\end{align*}
Using (\ref{fge1}) and Young inequality, we get
\begin{align*}
 I_{21}\leq C&\sum_{i,m=1}^3\sum_{j\neq i}\|\partial_i^2\partial_{i'} u_{i'}\|^\frac{1}{2} \left(\|\partial_i^2\partial_{i'} u_{i'}\|+\|\partial_{1i'}\partial_i^2 u_{i'}\|\right) ^{\frac{1}{2}} \|\partial_i^m u_j\|^\frac{1}{2} \left(\|\partial_i^m u_j\|+\|\partial_2\partial_i^m u_j\|\right)^\frac{1}{2} \|\partial_j\partial_i^{3-m} u_i\|^\frac{1}{2} \|\partial_{3j}\partial_i^{3-m} u_i\|^\frac{1}{2}\\
\leq&C\|u\|_{H^3} \sum_{i=1}^3\left(\|\partial_i u_{i'}\|_{H^3}^2 +\|\partial_i u_{i''}\|_{H^3}^2\right).
\end{align*}
The last term $I_{22}$ can be estimated similarly as $I_{21}$.
\begin{align*}
I_{22}\leq&C\|u\|_{H^3} \sum_{i=1}^3\left(\|\partial_i u_{i'}\|_{H^3}^2 +\|\partial_i u_{i''}\|_{H^3}^2\right).
\end{align*}
For the term $I_3$, $k$ can choose $i'$ or $i''$. We do the following decomposition
\begin{align*}
I_3=&-\sum_{m=1}^2\sum_{i=1}^3\sum_{k\neq i} C_3^m\int \partial_i^m u_i \partial_i\partial_i^{3-m} u_k\partial_i^3 u_k dx -\sum_{i=1}^3\sum_{k\neq i} C_3^m\int \partial_i^m u_i \partial_i\partial_i^{3-m} u_k\partial_i^3 u_k dx:=I_{31}+I_{32}.
\end{align*}
For the term $I_{31}$, from (\ref{fge1}) and Young inequality, we obtain
\begin{align*}
I_{31}\leq&C\sum_{m=1}^2\sum_{i=1}^3\sum_{k\neq i}\|\partial_i^3 u_k\|^\frac{1}{2} \left(\|\partial_i^3 u_k\|+\|\partial_1 \partial_i^3 u_k\|\right)^\frac{1}{2}  \|\partial_i^{4-m} u_k\|^\frac{1}{2}  \left(\|\partial_i^{4-m} u_k\|^\frac{1}{2}+ \|\partial_2\partial_i^{4-m} u_k\|\right)^\frac{1}{2} \|\partial_i^m u_i\|^\frac{1}{2} \|\partial_3\partial_i^m u_i\|^\frac{1}{2}\\
\leq&C\|u\|_{H^3} \sum_{i=1}^3\left(\|\partial_i u_{i'}\|_{H^3}^2 +\|\partial_i u_{i''}\|_{H^3}^2\right).
\end{align*}
To estimate the term $I_{32}$, we can refer to $I_{31}$. Namely,
\begin{align*}
I_{32}\leq C\|u\|_{H^3} \sum_{i=1}^3\left(\|\partial_i u_{i'}\|_{H^3}^2 +\|\partial_i u_{i''}\|_{H^3}^2\right).
\end{align*}
Let $j=i'$,~$ k=i''$  or $j=i''$,~$ k=i'$. We use again inequality (\ref{fge1}) and Young inequality to get
\begin{align*}
I_4\leq&C\sum_{i,m=1}^3\sum_{j\neq i,k\neq i}\|\partial_i^3 u_k\|^\frac{1}{2} \left(\|\partial_i^3 u_k\|+\|\partial_1 \partial_i^3 u_k\|\right)^\frac{1}{2} \|\partial_i^m u_j\|^\frac{1}{2} \left(\|\partial_i^m u_j\|+ \|\partial_2\partial_i^m u_j\|\right)^\frac{1}{2} \|\partial_{j}\partial_i^{3-m} u_k\|^\frac{1}{2} \|\partial_{3j}\partial_i^{3-m} u_k\|^\frac{1}{2}\\
\leq&C\|u\|_{H^3} \sum_{i=1}^3\left(\|\partial_i u_{i'}\|_{H^3}^2 +\|\partial_i u_{i''}\|_{H^3}^2\right).
\end{align*}
Combining the above these estimates, we obtain
\begin{align}\label{I}
I\leq C\|u\|_{H^3} \sum_{i=1}^3\left(\|\partial_{i'} u_i\|_{H^3}^2+\|\partial_{i''} u_i\|_{H^3}^2\right).
\end{align}
\noindent\textbf {Estimate of $J$.} We write
\begin{align*}
J=\sum_{m=1}^3\sum_{i=1}^2 C_3^m\int\partial_i^m b\cdot\nabla \partial_{i}^{3-m}b\cdot\partial_i^3 u dx+ \sum_{m=1}^3C_3^m\int\partial_3^m b\cdot\nabla \partial_3^{3-m}b\cdot\partial_3^3 u dx:=J_1+J_2.
\end{align*}
For the first term $J_1$, we split it into two parts
\begin{align*}
J_1=&\sum_{m,i,k=1}^2 C_3^m \int\partial_i^m b\cdot\nabla \partial_{i}^{3-m}b_k\partial_i^3 u_k dx+ \sum_{m,i=1}^2 C_3^m \int\partial_i^m b\cdot\nabla \partial_{i}^{3-m}b_3\partial_i^3 u_3 dx\\&+\sum_{i,j=1}^2 \int\partial_i^3 b_j \partial_j b\partial_i^3 udx+\sum_{i=1}^2 \int\partial_i^3 b_3\partial_3 b\partial_i^3 udx:=J_{11}+J_{12}+J_{13}+J_{14}.
\end{align*}
To estimate the term $J_{11}$, taking advantage of (\ref{fge1}) and Young inequality, we get
\begin{align*}
J_{11}\leq&C\sum_{m,i,k=1}^2\|\partial_i^m b\|^\frac{1}{2} \left(\|\partial_i^m b\|+\|\partial_1 \partial_i^m b\|\right)^\frac{1}{2} \|\nabla \partial_{i}^{3-m}b_k\|^\frac{1}{2} \left(\|\nabla \partial_{i}^{3-m}b_k\|^\frac{1}{2} +\|\partial_2\nabla \partial_{i}^{3-m}b_k\|\right)^\frac{1}{2} \|\partial_i^3 u_k\|^\frac{1}{2} \|\partial_3\partial_i^3 u_k\|^\frac{1}{2}\\
\leq&C\|(u,b)\|_{H^3} \left(\|\partial_3 u_1\|_{H^3}^2+\|\partial_3 u_2\|_{H^3}^2+\|\partial_1 b\|_{H^3}^2+\|\partial_2 b\|_{H^3}^2\right).
\end{align*}
In a similar way, this yields
\begin{align*}
J_{12}\leq&C\|(u,b)\|_{H^3} \left(\|\partial_3 u_1\|_{H^3}^2+\|\partial_3 u_2\|_{H^3}^2+\|\partial_1 b\|_{H^3}^2+\|\partial_2 b\|_{H^3}^2\right).
\end{align*}
To estimate the $J_{13}$, we decompose it into two terms
\begin{align*}
J_{13}=&\sum_{i,j,k=1}^2 \int\partial_i^3 b_j\partial_jb_k\partial_i^3 u_k dx+ \sum_{i,j=1}^2 \int\partial_i^3 b_j\partial_jb_3\partial_i^3 u_3 dx:=J_{131}+J_{132}.
\end{align*}
Using again (\ref{fge1}) and Young inequality, we can infer
\begin{align*}
J_{131}\leq&C\sum_{i,j,k=1}^2\|\partial_i^3 b_j\|^\frac{1}{2} \left(\|\partial_i^3 b_j\|+\|\partial_1 \partial_i^3 b_j\|\right)^\frac{1}{2} \|\partial_j b_k\|^\frac{1}{2} \left(\|\partial_j b_k\|^\frac{1}{2} +\|\partial_2\partial_j b_k\|\right)^\frac{1}{2} \|\partial_i^3 u_k\|^\frac{1}{2} \|\partial_3\partial_i^3 u_k\|^\frac{1}{2}\\
\leq&C\|(u,b)\|_{H^3} \sum_{i,j=1}^2\left(\|\partial_i u_3\|_{H^3}^2+\|\partial_3 u_i\|_{H^3}^2+\|\partial_i b_j\|_{H^3}^2\right).
\end{align*}
Similarly,
\begin{align*}
J_{132}\leq C\|(u,b)\|_{H^3} \sum_{i,j=1}^2\left(\|\partial_i u_3\|_{H^3}^2+\|\partial_3 u_i\|_{H^3}^2+\|\partial_i b_j\|_{H^3}^2\right).
\end{align*}
To estimate the $J_{14}$, we decompose
\begin{align*}
J_{14}=\sum_{i=1}^2 \int\partial_i^3 b_3\partial_3 b_i\partial_i^3 u_i dx+\int\partial_1^3 b_3\partial_3 b_2\partial_1^3 u_2dx+\int\partial_2^3 b_3\partial_3 b_1\partial_2^3 u_1 dx :=J_{141}+J_{142}+J_{143}.
\end{align*}
By using $div ~u=0$ and (\ref{fge1}) and Young inequality, one has
\begin{align*}
J_{141}=&-\sum_{i=1}^2 \int\partial_i^3 b_3\partial_3 b_i\partial_i^2 (\partial_{i'} u_{i'}+\partial_{i''} u_{i''}) dx\\
\leq&C\sum_{i=1}^2\|\partial_i^3 b_3\|^\frac{1}{2} \left(\|\partial_i^3 b_3\|+\|\partial_1 \partial_i^3 b_3\|\right)^\frac{1}{2} \|\partial_3 b_i\|^\frac{1}{2} \left(\|\partial_3 b_i\|^\frac{1}{2} +\|\partial_{23} b_i\|\right)^\frac{1}{2} \\ &\times\left(\|\partial_i^2\partial_{i'} u_{i'}\|^\frac{1}{2} \|\partial_3\partial_i^2\partial_{i'} u_{i'}\|^\frac{1}{2}+\|\partial_i^2\partial_{i''} u_{i''}\|^\frac{1}{2} \|\partial_3\partial_i^2\partial_{i''} u_{i''}\|^\frac{1}{2}\right)\\
\leq&C\|b\|_{H^3} \sum_{i=1}^2\left(\|\partial_i u_{i'}\|_{H^3}^2+\|\partial_i u_{i''}\|_{H^3}^2+\|\partial_i b_3\|_{H^3}^2\right).
\end{align*}
The term $J_{142}$, $J_{143}$ can be treated by the way of the term $J_{131}$, we have
\begin{align*}
J_{142}+J_{143}\leq C\|b\|_{H^3} \left(\|\partial_1 u_2\|_{H^3}^2+\|\partial_2 u_1\|_{H^3}^2+\|\partial_1 b_3+\|\partial_2 b_3\|_{H^3}^2\right).
\end{align*}
For the term $J_2$, We do the decomposition
\begin{align*}
J_2=&\sum_{m=1}^3\sum_{j,k=1}^2 C_3^m \int\partial_3^m b_j\partial_j\partial_3^{3-m}b_k \partial_3^3 u_k dx +\sum_{m=1}^3\sum_{j=1}^2 C_3^m \int\partial_3^m b_j\partial_j\partial_3^{3-m}b_3 \partial_3^3 u_3 dx\\ &+\sum_{m=1}^3\sum_{k=1}^2 C_3^m \int\partial_3^m b_3\partial_3^{4-m}b_k \partial_3^3 u_k dx +\sum_{m=1}^3 C_3^m \int\partial_3^m b_3\partial_3^{4-m}b_3 \partial_3^3 u_3 dx:=J_{21}+J_{22}+J_{23}+J_{24}.
\end{align*}
To bounded the first term $J_{21}$, we use (\ref{fge1}) and Young inequality to get
\begin{align*}
J_{21}\leq&C\sum_{m=1}^3\sum_{j,k=1}^2\|\partial_3^m b\|^\frac{1}{2}\left(\|\partial_3^m b\|^\frac{1}{2} +\|\partial_1\partial_3^m b\|\right)^\frac{1}{2} \|\partial_j \partial_3^{3-m}b_3\|^\frac{1}{2}\left(\|\partial_j \partial_3^{3-m}b_3\| +\|\partial_{2j} \partial_3^{3-m}b_3\|\right)^\frac{1}{2}\|\partial_3^3 u_k\|^\frac{1}{2} \|\partial_3^4 u_k\|^\frac{1}{2} \\
\leq&C\|b\|_{H^3} \sum_{k=1}^2\left(\|\partial_3 u_k\|_{H^3}^2+\|\partial_1 b_k\|_{H^3}^2+\|\partial_2 b_k\|_{H^3}^2\right).
\end{align*}
Similarly,
\begin{align*}
J_{22}\leq C\|b\|_{H^3} \sum_{k=1}^2\left(\|\partial_3 u_k\|_{H^3}^2+\|\partial_1 b\|_{H^3}^2+\|\partial_2 b_3\|_{H^3}^2\right).
\end{align*}
By using the divergence free condition $div~b=0$, inequality (\ref{fge1}) and Young inequality, we obtain
\begin{align*}
J_{23}=&-\sum_{m=1}^3\sum_{i,k=1}^2 C_3^m\int\partial_3^{m-1}\partial_i b_i \partial_3^{4-m}b_k \partial_3^3 u_k dx\\
\leq&C\sum_{m=1}^3\sum_{i,k=1}^2\|\partial_3^{m-1} \partial_i b_i\|^\frac{1}{2} \left(\|\partial_3^{m-1} \partial_i b_i\|+\|\partial_{1i}\partial_3^{m-1}  b_i\|\right)^\frac{1}{2} \|\partial_3^{4-m}b_k\|^\frac{1}{2} \left(\|\partial_3^{4-m}b_k\|+\|\partial_2\partial_3^{4-m}b_k\|\right)^\frac{1}{2} \|\partial_3^3 u_k\|^\frac{1}{2} \|\partial_3^4 u_k\|^\frac{1}{2}\\
\leq&C\|b\|_{H^3} \sum_{i=1}^2\left(\|\partial_3 u_i\|_{H^3}^2+\|\partial_1 b_i\|_{H^3}^2+\|\partial_2 b_i\|_{H^3}^2\right).
\end{align*}
In the similarly way,
\begin{align*}
J_{24}\leq&C\|b\|_{H^3} \sum_{i=1}^2\left(\|\partial_3 u_i\|_{H^3}^2+\|\partial_i b_i\|_{H^3}^2+\|\partial_2 b_3\|_{H^3}^2\right).
\end{align*}
Combining these estimates,
\begin{align}\label{J}
J\leq C\|b\|_{H^3} \sum_{i=1}^3\left(\|\partial_{i'} u_i\|_{H^3}^2+\|\partial_{i''} u_i\|_{H^3}^2+\|\partial_1 b\|_{H^3}^2+\|\partial_2 b\|_{H^3}^2\right).
\end{align}
\noindent\textbf {Estimate of $K$.} We write
\begin{align*}
K=&-\sum_{m=1}^3\sum_{i,j=1}^2 C_3^m\int \partial_i^m u_j \partial_j\partial_i^{3-m} b\cdot\partial_i^3 b dx-\sum_{m=1}^3\sum_{i=1}^2 C_3^m\int \partial_i^m u_3 \partial_3\partial_i^{3-m} b\cdot\partial_i^3 b dx\\&-\sum_{m=1}^3 C_3^m\int \partial_3^m u\cdot\nabla \partial_3^{3-m} b\partial_3^3 b dx:=K_1+K_2+K_3.
\end{align*}
For the first term $K_1$, by using (\ref{fge1}) and Young inequality, we find that
\begin{align*}
K_1\leq&C\sum_{m=1}^3\sum_{i,j=1}^2\|\partial_i^3 b\|^\frac{1}{2} \left(\|\partial_i^3 b\|+\|\partial_1\partial_i^3 b\|\right)^\frac{1}{2} \|\partial_j\partial_i^{3-m} b\|^\frac{1}{2} \left(\|\partial_j\partial_i^{3-m} b\|+\|\partial_{2j}\partial_i^{3-m} b\|\right)^\frac{1}{2} \|\partial_i^m u_j\|^\frac{1}{2} \|\partial_3\partial_i^m u_j\|^\frac{1}{2}\\
\leq&C\|(u,b)\|_{H^3} \sum_{i=1}^2\left(\|\partial_3 u_i\|_{H^3}^2+\|\partial_{i} b\|_{H^3}^2\right).
\end{align*}
Similarly,
\begin{align*}
K_2\leq C\|(u,b)\|_{H^3} \sum_{i=1}^2\left(\|\partial_i u_3\|_{H^3}^2+\|\partial_{i} b\|_{H^3}^2\right).
\end{align*}
For the last term $K_3$, we decompose it into four parts
\begin{align*}
K_3=&-\sum_{m=1}^3 C_3^m\int \partial_3^m u_3\partial_3 \partial_3^{3-m} b_3\partial_3^3 b_3 dx-\sum_{m=1}^3\sum_{j=1}^2 C_3^m\int \partial_3^m u_j\partial_j \partial_3^{3-m} b_3\partial_3^3 b_3dx\\
&-\sum_{m=1}^3\sum_{k=1}^2 C_3^m\int \partial_3^m u_3\partial_3 \partial_3^{3-m} b_k\partial_3^3 b_k dx-\sum_{m=1}^3\sum_{j,k=1}^2 C_3^m\int \partial_3^m u_j\partial_j \partial_3^{3-m} b_k\partial_3^3 b_kdx\\:=&K_{31}+K_{32}+K_{33}+K_{34}.
\end{align*}
Thanks to the divergence free condition $div~b=0$, (\ref{fge1}) and Young inequality, we get
\begin{align*}
K_{31}=&-\sum_{m=1}^3\sum_{j,k=1}^2 C_3^m\int \partial_3^m u_3 \partial_3^{3-m} \partial_jb_j\partial_3^2\partial_k b_k dx\\
\leq&C\sum_{m=1}^3\sum_{j,k=1}^2\|\partial_3^m u_3\|^\frac{1}{2} \left(\|\partial_3^m u_3\|+\|\partial_1\partial_3^m u_3\|\right)^\frac{1}{2} \|\partial_3^2\partial_k b_k\|^\frac{1}{2} \left(\|\partial_3^2\partial_k b_k\|+\|\partial_{2k}\partial_3^2 b_k\|\right)^\frac{1}{2} \|\partial_3^{3-m} \partial_jb_j\|^\frac{1}{2} \|\partial_{3j}\partial_3^{3-m} b_j\|^\frac{1}{2}\\
\leq&C\|(u,b)\|_{H^3} \sum_{k=1}^2\left(\|\partial_1 u_3\|_{H^3}^2+\|\partial_k b_k\|_{H^3}^2\right).
\end{align*}
Similarly,
\begin{align*}
K_{32}\leq&C\|(u,b)\|_{H^3} \sum_{k=1}^2\left(\|\partial_3 u_k\|_{H^3}^2+\|\partial_k b\|_{H^3}^2\right).
\end{align*}
To estimate the term $K_{33}$, we decompose
\begin{align*}
K_{33}=&\sum_{m=1}^3\sum_{i,k=1}^2 C_3^m\int \partial_3^{m-1}\partial_i u_i\partial_3^{4-m} b_k\partial_3^3 b_k dx=\sum_{m=1}^3\sum_{i,k=1}^2 C_3^m\int \partial_3^{m-1}\partial_i (\overline{u}_i+\widetilde{u}_i)\partial_3^{4-m} (\overline{b}_k+\widetilde{b}_k)\partial_3^3 (\overline{b}_k+\widetilde{b}_k) dx\\
=&\sum_{m=1}^3\sum_{i,k=1}^2 C_3^m\int \partial_3^{m-1}\partial_i \widetilde{u}_i\partial_3^{4-m} \overline{b}_k\partial_3^3 \widetilde{b}_k dx+\sum_{m=1}^3\sum_{i,k=1}^2 C_3^m\int \partial_3^{m-1}\partial_i \widetilde{u}_i\partial_3^{4-m} \widetilde{b}_k\partial_3^3 \overline{b}_k dx\\&+\sum_{m=1}^3\sum_{i,k=1}^2 C_3^m\int \partial_3^{m-1}\partial_i \widetilde{u}_i\partial_3^{4-m} \widetilde{b}_k\partial_3^3 \widetilde{b}_k dx:=K_{331}+K_{332}+K_{333}.
\end{align*}
For the term $K_{331}$, we split it into three parts
\begin{align*}
K_{331}=&3\sum_{k=1}^2 \int \partial_1 \widetilde{u}_1\partial_3^3 \overline{b}_k\partial_3^3 \widetilde{b}_k dx+3\sum_{k=1}^2 \int \partial_2 \widetilde{u}_2\partial_3^3 \overline{b}_k\partial_3^3 \widetilde{b}_k dx\\ &+\sum_{m=2}^3\sum_{i,k=1}^2 C_3^m\int \partial_3^{m-1}\partial_i \widetilde{u}_i\partial_3^{4-m} \overline{b}_k\partial_3^3 \widetilde{b}_k dx:=K_{3311}+K_{3312}+K_{3313}.
\end{align*}
It follows by using the inequality (\ref{fge2}), Lemma (\ref{lempoincare}) and Young inequality that
\begin{align*}
K_{3311}\leq&C\sum_{k=1}^2\|\partial_3^3 \widetilde{b}_k\|^\frac{1}{2} \|\partial_1\partial_3^3 \widetilde{b}_k\|^\frac{1}{2} \|\partial_3^3 \overline{b}_k\|^\frac{1}{2} \left(\|\partial_3^3 \overline{b}_k\| +\|\partial_2\partial_3^3 \overline{b}_k\|\right)^\frac{1}{2} \|\partial_1\widetilde{u}_1\|^\frac{1}{2} \|\partial_{31}\widetilde{u}_1\|^\frac{1}{2}\\
\leq&C\sum_{k=1}^2 \|\partial_1\partial_3^3 \widetilde{b}_k\| \|\partial_3^3 \overline{b}_k\|^\frac{1}{2} \left(\|\partial_3^3 \overline{b}_k\| +\|\partial_2\partial_3^3 \overline{b}_k\|\right)^\frac{1}{2} \|\partial_{31}\widetilde{u}_1\|\\
\leq&C\|(u,b)\|_{H^3} \sum_{k=1}^2\left(\|\partial_2 u_1\|_{H^3}^2+\|\partial_3 u_1\|_{H^3}^2+\|\partial_1 b_k\|_{H^3}^2+\|\partial_2 b_k\|_{H^3}^2\right).
\end{align*}
The terms $K_{3311}$ and $K_{3312}$,~$K_{3313}$ can be estimated in same way,
\begin{align*}
K_{3312}\leq C\|(u,b)\|_{H^3} \sum_{k=1}^2\left(\|\partial_1 u_2\|_{H^3}^2+\|\partial_3 u_2\|_{H^3}^2+\|\partial_1 b_k\|_{H^3}^2+\|\partial_2 b_k\|_{H^3}^2\right)
\end{align*}
and
\begin{align*}
K_{3313}\leq C\|(u,b)\|_{H^3} \sum_{k=1}^2\left(\|\partial_3 u_k\|_{H^3}^2+\|\partial_1 b_k\|_{H^3}^2+\|\partial_2 b_k\|_{H^3}^2\right).
\end{align*}
For the term $K_{332}$, we rewrite
\begin{align*}
K_{332}=&3\sum_{i,k=1}^2 \int \partial_i \widetilde{u}_i\partial_3^3 \widetilde{b}_k \partial_3^3 \overline{b}_k dx+\sum_{m=2}^3\sum_{i,k=1}^2 C_3^m\int \partial_3^{m-1}\partial_i \widetilde{u}_i\partial_3^{4-m} \widetilde{b}_k\partial_3^3 \overline{b}_k dx:=K_{3321}+K_{3322}.
\end{align*}
The terms $K_{3321}$ and $K_{3322}$ can be bounded in same way to $K_{3311}$ and $K_{3313}$, respectively. we obtain
\begin{align*}
K_{3321}\leq C\|(u,b)\|_{H^3} \sum_{k=1}^2\left(\|\partial_2 u_1\|_{H^3}^2+\|\partial_3 u_1\|_{H^3}^2+\|\partial_1 u_2\|_{H^3}^2+\|\partial_3 u_2\|_{H^3}^2+\|\partial_1 b_k\|_{H^3}^2+\|\partial_2 b_k\|_{H^3}^2\right)
\end{align*}
and
\begin{align*}
K_{3313}\leq C\|(u,b)\|_{H^3} \sum_{k=1}^2\left(\|\partial_3 u_k\|_{H^3}^2+\|\partial_1 b_k\|_{H^3}^2+\|\partial_2 b_k\|_{H^3}^2\right).
\end{align*}
From the inequality (\ref{fge2}), Lemma (\ref{lempoincare}) and Young inequality, we can get
\begin{align*}
K_{333}=&3\sum_{i,k=1}^2 \int \partial_i \widetilde{u}_i\partial_3^3 \widetilde{b}_k \partial_3^3 \widetilde{b}_k dx+\sum_{m=2}^3\sum_{i,k=1}^2 C_3^m\int \partial_3^{m-1}\partial_i \widetilde{u}_i\partial_3^{4-m} \widetilde{b}_k\partial_3^3 \widetilde{b}_k dx\\
\leq&C\sum_{m=1}^3\sum_{i,k=1}^2\|\partial_3^{4-m} \widetilde{b}_k\|^\frac{1}{2} \|\partial_1\partial_3^{4-m} \widetilde{b}_k\|^\frac{1}{2} \|\partial_3^3 \widetilde{b}_k\|^\frac{1}{2} \|\partial_2\partial_3^3 \widetilde{b}_k\|^\frac{1}{2} \|\partial_3^{m-1}\partial_i \widetilde{u}_i\|^\frac{1}{2} \|\partial_3^m\partial_i \widetilde{u}_i\|^\frac{1}{2}\\
\leq&C\|(u,b)\|_{H^3} \sum_{k=1}^2\left(\|\partial_3 u_k\|_{H^3}^2+\|\partial_1 b_k\|_{H^3}^2+\|\partial_2 b_k\|_{H^3}^2\right).
\end{align*}
Combining these cases, we obtain the bounded of the $K$.
Namely,
\begin{align}\label{K}
K\leq C\|(u,b)\|_{H^3} \sum_{i=1}^3\left(\|\partial_{i'} u_i\|_{H^3}^2+\|\partial_{i''} u_i\|_{H^3}^2+\|\partial_1 b\|_{H^3}^2+\|\partial_2 b\|_{H^3}^2\right).
\end{align}
\noindent\textbf {Estimate of $L$.} We split it into four parts
\begin{align*}
L=&\sum_{m=1}^3\sum_{i,k=1}^2  C_3^m\int \partial_i^m b\cdot\nabla  \partial_i^{3-m}u_k\partial_i^3 b_k dx+\sum_{m=1}^3\sum_{i=1}^2C_3^m\int \partial_i^m b\cdot\nabla  \partial_i^{3-m}u_3\partial_i^3 b_3 dx\\&+\sum_{m=1}^3\sum_{k=1}^2  C_3^m\int \partial_3^m b\cdot\nabla  \partial_3^{3-m}u_k\partial_i^3 b_k dx+\sum_{m=1}^3  C_3^m\int \partial_3^m b\cdot\nabla  \partial_3^{3-m}u_3\partial_i^3 b_3 dx:=L_1+L_2+L_3+L_4.
\end{align*}
The term $L_1$ and $L_2$ can be bounded in same way to $K_1$, we get
\begin{align*}
L_1\leq&C\|(u,b)\|_{H^3} \left(\|\partial_3 u_1\|_{H^3}^2+\|\partial_3 u_2\|_{H^3}^2+\|\partial_1 b\|_{H^3}^2+\|\partial_2 b\|_{H^3}^2\right)
\end{align*}
and
\begin{align*}
L_2\leq&C\|(u,b)\|_{H^3} \left(\|\partial_2 u_3\|_{H^3}^2+\|\partial_1 b\|_{H^3}^2+\|\partial_2 b\|_{H^3}^2\right).
\end{align*}
For the term $L_3$, we write
\begin{align*}
L_3=&\sum_{m=1}^3\sum_{j,k=1}^2 C_3^m\int \partial_3^m b_j\partial_j \partial_3^{3-m}u_k\partial_3^3 b_k dx+\sum_{m=1}^3\sum_{k=1}^2 C_3^m\int \partial_3^m b_3 \partial_3^{4-m}u_k\partial_3^3 b_k dx:=L_{31}+L_{32}.
\end{align*}
For the term $L_{31}$, we do the following decomposition 
\begin{align*}
L_{31}=&\sum_{m=1}^3\sum_{j,k=1}^2 C_3^m\int \partial_3^m \overline{b}_j\partial_j \partial_3^{3-m}\widetilde{u}_k\partial_3^3 \widetilde{b}_k dx +\sum_{m=1}^3\sum_{j,k=1}^2 C_3^m\int \partial_3^m\widetilde{b}_j\partial_j \partial_3^{3-m}\widetilde{u}_k\partial_3^3 \overline{b}_k dx\\ &+\sum_{m=1}^3\sum_{j,k=1}^2 C_3^m\int \partial_3^m\widetilde{b}_j\partial_j \partial_3^{3-m}\widetilde{u}_k\partial_3^3 \widetilde{b}_k dx
:= L_{311}+L_{312}+L_{313}.
\end{align*}
To bounded the term $L_{311}$, we rewrite
\begin{align*}
L_{311}=&\sum_{m,j,k=1}^2 C_3^m\int \partial_3^m \overline{b}_j\partial_j \partial_3^{3-m}\widetilde{u}_k\partial_3^3 \widetilde{b}_k dx +\sum_{j,k=1}^2 \int \partial_3^3 \overline{b}_j\partial_j \widetilde{u}_k\partial_3^3 \widetilde{b}_k dx:=L_{3111}+L_{3112}.
\end{align*}
The term $L_{3111}$ can be treated by the way of the term $K_{3313}$. we obtain
\begin{align*}
L_{3111}\leq&C\|b\|_{H^3} \left(\|\partial_3 u_1\|_{H^3}^2+\|\partial_3 u_2\|_{H^3}^2+\|\partial_2 b_1\|_{H^3}^2+\|\partial_2 b_2\|_{H^3}^2\right).
\end{align*}
For the term $L_{3112}$, we split it into four part
\begin{align*}
L_{3112}=& \int \partial_3^3 \overline{b}_1\partial_1 \widetilde{u}_1\partial_3^3 \widetilde{b}_1 dx+\int \partial_3^3 \overline{b}_2\partial_2 \widetilde{u}_2\partial_3^3 \widetilde{b}_2 dx+\int \partial_3^3 \overline{b}_1\partial_1 \widetilde{u}_2\partial_3^3 \widetilde{b}_2 dx\\
&+\int \partial_3^3 \overline{b}_2\partial_2 \widetilde{u}_1\partial_3^3 \widetilde{b}_1 dx
:= L_{31121}+L_{31122}+L_{31123}+L_{31124}.
\end{align*}
Using the inequality (\ref{fge2}), Lemma (\ref{lempoincare}) and Young inequality, we have
\begin{align*}
L_{31121}\leq&C \|\partial_3^3 \overline{b}_1\|^\frac{1}{2} \left(\|\partial_3^3 \overline{b}_1\|^\frac{1}{2} +\|\partial_1\partial_3^3 \overline{b}_1\|\right)^\frac{1}{2} \|\partial_3^3 \widetilde{b}_1\|^\frac{1}{2} \|\partial_2\partial_3^3 \widetilde{b}_1\|^\frac{1}{2} \|\partial_1\widetilde{u}_1\|^\frac{1}{2} \|\partial_{31}\widetilde{u}_1\|^\frac{1}{2}\\
\leq&C\|b\|_{H^3} \left(\|\partial_2 u_1\|_{H^3}^2+\|\partial_3 u_1\|_{H^3}^2+\|\partial_1 b_1\|_{H^3}^2+\|\partial_2 b_1\|_{H^3}^2\right).
\end{align*}
Similarly,
\begin{align*}
L_{31122}\leq& C\|b\|_{H^3} \left(\|\partial_1 u_2\|_{H^3}^2+\|\partial_3 u_2\|_{H^3}^2+\|\partial_1 b_2\|_{H^3}^2+\|\partial_2 b_2\|_{H^3}^2\right),\\
L_{31123}\leq& C\|b\|_{H^3} \left(\|\partial_1 u_2\|_{H^3}^2+\|\partial_1 b_1\|_{H^3}^2+\|\partial_2 b_2\|_{H^3}^2\right),\\
L_{31124}\leq& C\|b\|_{H^3} \left(\|\partial_2 u_1\|_{H^3}^2+\|\partial_1 b_1\|_{H^3}^2+\|\partial_2 b_2\|_{H^3}^2\right).
\end{align*}
The estimates for the term $L_{312}$ and $L_{313}$ are similar to that for the term $L_{311}$ and $K_{333}$, respectively.
To estimate the last term $L_{32}$, we use the same argument as the term $K_{31}$. One has
\begin{align*}
L_{32}\leq&C\|b\|_{H^3}  \sum_{k=1}^2\left(\|\partial_3 u_k\|_{H^3}^2+\|\partial_1 b_k\|_{H^3}^2+\|\partial_2 b_2\|_{H^3}^2\right).
\end{align*}
By using the divergence free condition $div~u=div~b=0$, we decompose
\begin{align*}
L_4=-\sum_{m=1}^3\sum_{j,k=1}^2  C_3^m\int \partial_3^m b_j\partial_j  \partial_3^{3-m}u_3\partial_i^2 \partial_kb_k dx+\sum_{m=1}^3\sum_{i,k=1}^2  C_3^m\int \partial_3^{m-1}\partial_i b_i  \partial_3^{4-m}u_3\partial_3^2 \partial_kb_k dx:=L_{41}+L_{42}.
\end{align*}
For the estimates of the term $L_{41}$ and $L_{42}$, we use the same idea as $K_1$. We get
\begin{align*}
L_{41}\leq&C\|b\|_{H^3} \sum_{k=1}^2\left(\|\partial_k u_3\|_{H^3}^2+\|\partial_1 b_k\|_{H^3}^2+\|\partial_2 b_2\|_{H^3}^2\right),\\
L_{42}\leq&C\|(u,b)\|_{H^3}\sum_{i=1}^2\left(\|\partial_3 u_i\|_{H^3}^2+\|\partial_1 b_i\|_{H^3}^2+\|\partial_2 b_i\|_{H^3}^2\right).
\end{align*}
Combining these estimates, we find that
\begin{align}\label{L}
L\leq C\|(u,b)\|_{H^3} \sum_{i=1}^3\left(\|\partial_{i'} u_i\|_{H^3}^2+\|\partial_{i''} u_i\|_{H^3}^2+\|\partial_1 b\|_{H^3}^2+\|\partial_2 b\|_{H^3}^2\right).
\end{align}
Gathering these estimates (\ref{I}), (\ref{J}), (\ref{K})~and the last estimate (\ref{L}), we finally get
\begin{align*}
\frac{1}{2}&\frac{d}{d t}\|(u, b)\|_{H^3}^2+\nu\sum_{i=1}^3\|(\partial_{i'} ,\partial_{i''}) u_i\|_{H^3}^2+\mu\|(\partial_1, \partial_2) b\|_{H^3}^2\\
&\leq C\|(u,b)\|_{H^3} \left(\sum_{i=1}^3\|(\partial_{i'} ,\partial_{i''}) u_i\|_{H^3}^2+\|(\partial_1,\partial_2) b\|_{H^3}^2\right).
\end{align*}
Integrating in time we obtain
\begin{align*}
&\sup_{\tau\in [0, t]}\|(u, b)\|_{H^3}^2+\nu\sum_{i=1}^3\int_0^t\|(\partial_{i'} ,\partial_{i''}) u_i\|_{H^3}^2 d\tau+\mu\int_0^t\|(\partial_1, \partial_2) b\|_{H^3}^2 d\tau\\
&\leq C\|(u_{in}, b_{in})\|_{H^3}^2+ C\sup_{\tau\in [0, t]}\|(u,b)\|_{H^3} \int_0^t\left(\sum_{i=1}^3\|(\partial_{i'} ,\partial_{i''}) u_i\|_{H^3}^2+\|(\partial_1,\partial_2) b\|_{H^3}^2\right)d\tau.
\end{align*}
We deduce that
\[
E(t)\leq E(0)+E(t)^{\frac{3}{2}}, \quad \forall ~ t\geq 0.
\]
We finish the proof for the global stability of the system (\ref{smhe}). Then we prove the uniqueness. Assume that $(u^{(1)}, p^{(1)}, b^{(1)})$ and $(u^{(2)}, p^{(2)}, b^{(2)})$ are two pairs of solutions of the system (\ref{smhe}) with the same initial data $(u_{in}, b_{in})$ on $[0, t]$. Denote
$$ u^{\delta}=u^{(1)}-u^{(2)}, ~p^{\delta}=p^{(1)}-p^{(2)}, ~b^{\delta}=b^{(1)}-b^{(2)}.$$
Each equation corresponds to a difference, we have
\begin{eqnarray}\label{smhedi}
\begin{cases}
\partial_{t} u^{\delta}+ u^{(1)}\cdot\nabla u^{\delta} + u^{\delta} \cdot\nabla u^{(2)}=\nu\begin{bmatrix} \partial_{22} + \partial_{33}\\ \partial_{11} + \partial_{33}\\ \partial_{11} + \partial_{22} \end{bmatrix}~u^{\delta}-\nabla p^{\delta}+b^{(1)} \cdot\nabla b^{\delta}+b^{\delta} \cdot\nabla b^{(2)} +\partial_3 b^{\delta}, \\
\partial_{t} b^{\delta}+ u^{(1)} \cdot\nabla b^{\delta}+u^{\delta} \cdot\nabla b^{(2)}~=\mu\left(\partial_{11} + \partial_{22} \right) b^{\delta} +b^{(1)}\cdot\nabla u^{\delta}+b^{\delta}\cdot\nabla u^{(2)} +\partial_3 u^{\delta}, \\
div~ u^{\delta}=div~ b^{\delta}=0, \\
(u^{\delta}, b^{\delta})\vert _{t=0}=(0, 0). \\
\end{cases}
\end{eqnarray}
Taking the $L^2$ inner product of the system (\ref{smhedi}) with $(u^{\delta}, b^{\delta})$, we get that
\begin{align*}
\frac{1}{2}\frac{d}{dt}&\|(u^{\delta}, b^{\delta})\|^2+\nu\left(\|(\partial_2, \partial_3) u^{\delta}_1\|^2
+\|(\partial_1, \partial_3) u^{\delta}_2\|^2+\|(\partial_1, \partial_2) u^{\delta}_3\|^2\right)+\mu\|(\partial_1, \partial_2) b^{\delta}\|^2\\
&=-\int u^{\delta}\cdot\nabla u^{(2)}\cdot u^{\delta} dx+\int b^{\delta}\cdot\nabla b^{(2)}\cdot u^{\delta} dx-\int u^{\delta}\cdot\nabla b^{(2)}\cdot b^{\delta} dx+\int b^{\delta}\cdot\nabla u^{(2)}\cdot b^{\delta} dx  \\&:=M_1+M_2+M_3+M_4.
\end{align*}
We divided $M_1$ into four terms
\begin{align*}
M_1=&-\sum_{k=1}^2\int u^{\delta}_1\partial_1 u^{(2)}_k u^{\delta}_k dx -\sum_{j=2}^3\sum_{k=1}^2\int u^{\delta}_j\partial_j u^{(2)}_k u^{\delta}_k dx -\sum_{j=1}^2\int u^{\delta}_j\partial_j u^{(2)}_3 u^{\delta}_3 dx\\& -\int u^{\delta}_3\partial_3 u^{(2)}_3 u^{\delta}_3 dx:=M_{11}+M_{12}+M_{13}+M_{14}.
\end{align*}
For the item $M_1$, by using inequality (\ref{fge1}) and Young inequality, we have
\begin{align*}
M_{11}\leq&C\sum_{k=1}^2\|\partial_1 u^{(2)}_k\|^\frac{1}{2} \left(\|\partial_1 u^{(2)}_k\|+\|\partial_1^2 u^{(2)}_k\|\right)^\frac{1}{2} \|u^{\delta}_1\|^\frac{1}{2} \left(\|u^{\delta}_1\|+\|\partial_2 u^{\delta}_1\|\right)^\frac{1}{2} \|u^{\delta}_k\|^\frac{1}{2} \|\partial_3 u^{\delta}_k\|^\frac{1}{2}\\
\leq&C\|u^{\delta}\|^2+\frac{\nu}{2}\|\partial_2 u^{\delta}_1\|^2+\frac{\nu}{10}\sum_{k=1}^2\|\partial_3 u^{\delta}_k\|^2.
\end{align*}
Similarly,
\begin{align*}
M_{12}\leq&C\|u^{\delta}\|^2+\frac{\nu}{2}\|\partial_1 u^{\delta}_2\|^2+\frac{\nu}{6}\|\partial_1 u^{\delta}_3\|^2+\frac{\nu}{10}\sum_{k=1}^2\|\partial_3 u^{\delta}_k\|^2,\\
M_{13}\leq&C\|u^{\delta}\|^2+\frac{\nu}{6}\|\partial_1 u^{\delta}_3\|^2+\frac{\nu}{10}\sum_{j=1}^2\|\partial_3 u^{\delta}_j\|^2,\\
M_{14}\leq&C\|u^{\delta}\|^2+\frac{\nu}{6}\|\partial_1 u^{\delta}_3\|^2+\frac{\nu}{6}\|\partial_2 u^{\delta}_3\|^2.
\end{align*}
To bounded the term $M_2$, We write that
\begin{align*}
M_2=&\sum_{k=1}^2\int b^{\delta}\cdot\nabla b^{(2)}_k u^{\delta}_k dx +\int b^{\delta}\cdot\nabla b^{(2)}_3 u^{\delta}_3 dx:=M_{21}+M_{22}.
\end{align*}
Thanks to inequality (\ref{fge1}) and Young inequality, we get
\begin{align*}
M_{21}\leq&C\sum_{k=1}^2\|b^{\delta}\|^\frac{1}{2} \left(\|b^{\delta}\|+\|\partial_1 b^{\delta}\|\right)^\frac{1}{2} \|\nabla b^{(2)}_k\|^\frac{1}{2} \left(\|\nabla b^{(2)}_k\|+\|\partial_2 \nabla b^{(2)}_k\|\right)^\frac{1}{2} \|u^{\delta}_k\|^\frac{1}{2} \|\partial_3 u^{\delta}_k\|^\frac{1}{2}\\
\leq&C\|(u^{\delta}, b^{\delta})\|^2+\frac{\nu}{10}\sum_{k=1}^2\|\partial_3 u^{\delta}_k\|^2+\frac{\mu}{10}\|\partial_1 b^{\delta}\|^2.
\end{align*}
and
\begin{align*}
M_{22}\leq&C\|(u^{\delta}, b^{\delta})\|^2+\frac{\nu}{6}\|\partial_2 u^{\delta}_3\|^2+\frac{\mu}{10}\|\partial_1 b^{\delta}\|^2.
\end{align*}
The term $M_3$, $M_4$ can be controlled by the same way.
One has
\begin{align*}
M_3\leq&C\|(u^{\delta}, b^{\delta})\|^2+\frac{\nu}{6}\|\partial_2 u^{\delta}_3\|^2+\frac{\nu}{10}\sum_{k=1}^2\|\partial_3 u^{\delta}_k\|^2+\frac{\mu}{5}\|\partial_1 b^{\delta}\|^2,\\
M_4\leq&C\|b^{\delta}\|^2+\frac{\mu}{10}\|\partial_1 b^{\delta}\|^2+\frac{\mu}{2}\|\partial_2 b^{\delta}\|^2.
\end{align*}
Combining this estimate, we get
\begin{align*}
\frac{1}{2}\frac{d}{dt}&\|(u^{\delta}, b^{\delta})\|^2+\nu\left(\|(\partial_2, \partial_3) u^{\delta}_1\|^2
+\|(\partial_1, \partial_3) u^{\delta}_2\|^2+\|(\partial_1, \partial_2) u^{\delta}_3\|^2\right)+\mu\|(\partial_1, \partial_2) b^{\delta}\|^2\\
\leq &C\|(u^{\delta}, b^{\delta})\|^2+\frac{\nu}{2}\left(\|(\partial_2, \partial_3) u^{\delta}_1\|^2
+\|(\partial_1, \partial_3) u^{\delta}_2\|^2+\|(\partial_1, \partial_2) u^{\delta}_3\|^2\right)+\frac{\mu}{2}\|(\partial_1, \partial_2) b^{\delta}\|^2.
\end{align*}
From Gr\"{o}nwall inequality, we get
\begin{align*}
\|(u^{\delta}, b^{\delta})\|^2=0.
\end{align*}
This finishes the proof of Theorem \ref{theorem1t}.
\end{proof}
\section{Proof of Theorem \ref{theorem2t}.}\label{sec4t}
We now verify the exponential decay of $\|(\widetilde{u}, \widetilde{b})\|_{H^2}$.
\begin{proof}
 Dotting (\ref{smhedv}) by $(\widetilde{u},\widetilde{b})$ and applying $\partial_i^2$ to (\ref{smhedv}) and taking $L^2$ inner product with $(\partial_i^2 \widetilde{u}, \partial_i^2 \widetilde{b})$. Then taking summation for $i$ from $1$ to $3$, we obtain
\begin{align*} 
\frac{1}{2}\frac{d}{dt}&\|(\widetilde{u},\widetilde{b})\|_{H^2}^2+\nu\left(\|\nabla \widetilde{u}_1\|_{H^2}^2
+\|\nabla \widetilde{u}_2\|_{H^2}^2+\|(\partial_1, \partial_2) \widetilde{u}_3\|_{H^2}^2\right)+\mu\|(\partial_1, \partial_2) \widetilde{b}\|_{H^2}^2\\
=&-\int \widetilde{u \cdot\nabla \widetilde{u}}\cdot \widetilde{u} dx-\int \widetilde{u}_3\partial_3 \overline{u}\cdot \widetilde{u} dx+\int \widetilde{b \cdot\nabla \widetilde{b}}\cdot \widetilde{u} dx+\int \widetilde{b}_3\partial_3 \overline{b}\cdot \widetilde{u} dx-\int \widetilde{u \cdot\nabla \widetilde{b}}\cdot \widetilde{b} dx-\int \widetilde{u}_3\partial_3 \overline{b}\cdot \widetilde{b} dx\\
&+\int \widetilde{b \cdot\nabla \widetilde{u}}\cdot \widetilde{b} dx+\int \widetilde{b}_3\partial_3 \overline{u}\cdot \widetilde{b} dx-\sum_{i=1}^3\int \partial_i^2\widetilde{u \cdot\nabla \widetilde{u}}\cdot \partial_i^2\widetilde{u} dx-\sum_{i=1}^3\int \partial_i^2(\widetilde{u}_3\partial_3 \overline{u})\cdot \partial_i^2\widetilde{u} dx\\
&+\sum_{i=1}^3\int \partial_i^2\widetilde{b \cdot\nabla \widetilde{b}}\cdot \partial_i^2\widetilde{u} dx+\sum_{i=1}^3\int \partial_i^2(\widetilde{b}_3\partial_3 \overline{b})\cdot \partial_i^2\widetilde{u} dx-\sum_{i=1}^3\int \partial_i^2\widetilde{u \cdot\nabla \widetilde{b}}\cdot \partial_i^2\widetilde{b} dx-\sum_{i=1}^3\int \partial_i^2(\widetilde{u}_3\partial_3 \overline{b})\cdot \partial_i^2\widetilde{b} dx\\
&+\sum_{i=1}^3\int \partial_i^2\widetilde{b \cdot\nabla \widetilde{u}}\cdot \partial_i^2\widetilde{b} dx+\sum_{i=1}^3\int \partial_i^2(\widetilde{b}_3\partial_3 \overline{u})\cdot \partial_i^2\widetilde{b} dx
:=\sum_{i=1}^8N_i+\sum_{i=1}^8Q_i.
\end{align*}
By Lemma (\ref{lempro}), we get
\begin{align*}
N_1&=-\int u \cdot\nabla \widetilde{u}\cdot \widetilde{u} dx+\int \overline{u \cdot\nabla \widetilde{u}}\cdot \widetilde{u} dx=-\int u \cdot\nabla \widetilde{u}\cdot \widetilde{u} dx=0,\\
N_5&=-\int u \cdot\nabla \widetilde{b}\cdot \widetilde{b} dx+\int \overline{u \cdot\nabla \widetilde{b}}\cdot \widetilde{b} dx=-\int u \cdot\nabla \widetilde{b}\cdot \widetilde{b} dx=0,\\
N_3&+N_7=\int b \cdot\nabla \widetilde{b}\cdot \widetilde{u} dx+\int b \cdot\nabla \widetilde{u}\cdot \widetilde{b} dx=0.
\end{align*}
For the term $N_2$, it is divided into two cases
\begin{align*}
N_2&=-\sum_{K=1}^2\int \widetilde{u}_3\partial_3 \overline{u}_k \widetilde{u}_k dx-\int \widetilde{u}_3\partial_3 \overline{u}_3 \widetilde{u}_3 dx:=N_{21}+N_{22}.
\end{align*}
Invoking the inequality (\ref{fge2}), Lemma (\ref{lempoincare}) and Young inequality, we have
\begin{align*}
N_{21}\leq&C\sum_{k=1}^2\|\widetilde{u}_3\|^\frac{1}{2} \|\partial_1 \widetilde{u}_3\|^\frac{1}{2} \|\partial_3 \overline{u}_k\|^\frac{1}{2} \left(\|\partial_3 \overline{u}_k\|+\|\partial_{23} \overline{u}_k\|\right)^\frac{1}{2} \|\widetilde{u}_k\|^\frac{1}{2} \|\partial_3 \widetilde{u}_k\|^\frac{1}{2}\\
\leq&C\|u\|_{H^2}\sum_{k=1}^2\left(\|\partial_1 \widetilde{u}_3\|^2+\|\partial_3 \widetilde{u}_k\|^2\right).
\end{align*}
From the inequality (\ref{fge3}), Lemma (\ref{lempoincare}) and Young inequality, we obtain
\begin{align*}
N_{22}\leq&C\|\widetilde{u}_3\|^\frac{1}{2} \|\partial_1 \widetilde{u}_3\|^\frac{1}{2} \|\widetilde{u}_3\|^\frac{1}{2} \|\partial_2 \widetilde{u}_3\|^\frac{1}{2}\|\partial_3\overline{u}_3\|^\frac{1}{2} \|\partial_3^2\overline{u}_3\|^\frac{1}{2}
\leq C\|u\|_{H^2}\left(\|\partial_1 \widetilde{u}_3\|^2+\|\partial_2 \widetilde{u}_3\|^2\right).
\end{align*}
The estimate of the term $N_4$ can be treated by the way of the term $N_2$. One has
\begin{align*}
N_4\leq C\|b\|_{H^2}\left(\|(\partial_2,\partial_3) \widetilde{u}_1\|^2+\|(\partial_1,\partial_3) \widetilde{u}_2\|^2+\|\partial_2 \widetilde{u}_3\|^2+\|\partial_1 \widetilde{b}_3\|^2\right).
\end{align*}
The term $N_6, N_8$ can be bounded similarly as $N_{22}$. We get
\begin{align*}
N_6\leq C\|b\|_{H^2}\left(\|\partial_1 \widetilde{u}_3\|^2+\|\partial_2 \widetilde{b}\|^2\right), ~N_8\leq C\|u\|_{H^2}\left(\|\partial_1 \widetilde{b}_3\|^2+\|\partial_2 \widetilde{b}\|^2\right).
\end{align*}
For the term $Q_1$, we rewrite that
\begin{align*}
Q_1=&-\sum_{i=1}^3\int \partial_i^2(u \cdot\nabla \widetilde{u})\cdot \partial_i^2\widetilde{u} dx\\
=&-\sum_{m=1}^2\sum_{i=1}^3 C_2^m\int \partial_i^m u \cdot\nabla \partial_i^{2-m}\widetilde{u}_1 \partial_i^2\widetilde{u}_1 dx-\sum_{m=1}^2\sum_{i=1}^3 C_2^m\int \partial_i^m u \cdot\nabla \partial_i^{2-m}\widetilde{u}_2 \partial_i^2\widetilde{u}_2 dx\\
&-\sum_{m=1}^2\sum_{i=1}^3 C_2^m\int \partial_i^m u \cdot\nabla \partial_i^{2-m}\widetilde{u}_3 \partial_i^2\widetilde{u}_3 dx
:=Q_{11}+Q_{12}+Q_{13}.
\end{align*}
The term $Q_1$ can be controlled by the way of the term $N_2$. We obtain
\begin{align*}
Q_1\leq C\|u\|_{H^3}\left(\|(\partial_2,\partial_3) \widetilde{u}_1\|_{H^2}^2+\|(\partial_1,\partial_3) \widetilde{u}_2\|_{H^2}^2+\|(\partial_1,\partial_2) \widetilde{u}_3\|_{H^2}^2\right).
\end{align*}
For the term $Q_2$, we decompose
\begin{align*}
Q_2=&-\sum_{m=1}^2\sum_{i=1}^3 C_2^m\int \partial_i^m \widetilde{u}_3 \partial_3^{3-m}\overline{u}_1 \partial_i^2\widetilde{u}_1 dx-\sum_{m=1}^2\sum_{i=1}^3\sum_{k=2}^3 C_2^m\int \partial_i^m \widetilde{u}_3 \partial_3^{3-m}\overline{u}_k \partial_i^2\widetilde{u}_k dx-\sum_{i=1}^3 \int \widetilde{u}_3 \partial_3^3\overline{u}_1 \partial_i^2\widetilde{u}_1 dx\\
&-\sum_{i=1}^3\sum_{k=2}^3 \int \widetilde{u}_3 \partial_3^3\overline{u}_k \partial_i^2\widetilde{u}_k dx
:=Q_{21}+Q_{22}+Q_{23}+Q_{24}.
\end{align*}
For $Q_{21}$ and $ Q_{22}$, we use the same way as $N_2$. This yields
\begin{align*}
Q_{21}+Q_{22}\leq C\|u\|_{H^3}\left(\|\partial_2 \widetilde{u}_1\|_{H^2}^2+\|\partial_1 \widetilde{u}_2\|_{H^2}^2+\|\partial_1 \widetilde{u}_3\|_{H^2}^2+\|\partial_2 \widetilde{u}_3\|_{H^2}^2\right).
\end{align*}
Thanks to the inequality (\ref{fge23}), Lemma (\ref{lempoincare})  and Young inequality, we obtain
\begin{align*}
Q_{23}\leq&C\sum_{i=1}^3\|\widetilde{u}_3\|^\frac{1}{4} \|\partial_1 \widetilde{u}_3\|^\frac{1}{4} \|\partial_3 \widetilde{u}_3\|^\frac{1}{4}\|\partial_{13} \widetilde{u}_3\|^\frac{1}{4}
\|\partial_i^2\widetilde{u}_1\|^\frac{1}{2} \|\partial_2\partial_i^2\widetilde{u}_1\|^\frac{1}{2}\|\partial_3^3 \overline{u}_1\|\\
\leq&C\|u\|_{H^3}\left(\|\partial_1 \widetilde{u}_3\|_{H^2}^2+\|\partial_2 \widetilde{u}_1\|_{H^2}^2\right).
\end{align*}
Similarly,
\begin{align*}
Q_{24}\leq C\|u\|_{H^3}\left(\|\partial_1 \widetilde{u}_2\|_{H^2}^2+\|\partial_1 \widetilde{u}_3\|_{H^2}^2+\|\partial_2 \widetilde{u}_3\|_{H^2}^2\right).
\end{align*}
For writing convenience, the term $Q_3$ and $Q_7$ can be estimated together.
\begin{align*}
Q_3+Q_7=&\sum_{i=1}^3\int \partial_i^2(b \cdot\nabla \widetilde{b})\cdot \partial_i^2\widetilde{u} dx+\sum_{i=1}^3\int \partial_i^2(b \cdot\nabla \widetilde{u})\cdot \partial_i^2\widetilde{b} dx\\ =&\sum_{m=1}^2\sum_{i=1}^3 C_2^m\int \partial_i^m b\cdot\nabla \partial_i^{2-m}\widetilde{b}_1 \partial_i^2\widetilde{u}_1 dx+\sum_{m=1}^2\sum_{i=1}^3 \sum_{k=2}^3C_2^m\int \partial_i^m b\cdot\nabla \partial_i^{2-m}\widetilde{b}_k \partial_i^2\widetilde{u}_k dx\\
&+\sum_{m=1}^2\sum_{i=1}^3 C_2^m\int \partial_i^m b\cdot\nabla \partial_i^{2-m}\widetilde{u}_1 \partial_i^2\widetilde{b}_1 dx+\sum_{m=1}^2\sum_{i=1}^3 \sum_{k=2}^3C_2^m\int \partial_i^m b\cdot\nabla \partial_i^{2-m}\widetilde{u}_k \partial_i^2\widetilde{b}_k dx
\\:=&Q_{31}+Q_{32}+Q_{71}+Q_{72}.
\end{align*}
By using the inequality (\ref{fge3}), Lemma (\ref{lempoincare}) and Young inequality, we obtain
\begin{align*}
Q_{31}\leq&C\|\partial_i^{2-m}\widetilde{b}_1\|^\frac{1}{2} \|\partial_1 \partial_i^{2-m}\widetilde{b}_1\|^\frac{1}{2} \|\partial_i^2\widetilde{u}_1\|^\frac{1}{2} \|\partial_2 \partial_i^2\widetilde{u}_1\|^\frac{1}{2}\|\partial_i^m b\|^\frac{1}{2} \|\partial_3\partial_i^m b\|^\frac{1}{2}\\
\leq&C\|b\|_{H^3}\left(\|\partial_2\widetilde{u}_1\|_{H^2}^2+\|\partial_1 \widetilde{b}\|_{H^2}^2\right).
\end{align*}
In the same way,
\begin{align*}
Q_{32}\leq& C\|b\|_{H^3}\left(\|\partial_1\widetilde{u}_2\|_{H^2}^2+\|\partial_1\widetilde{u}_3\|_{H^2}^2+\|\partial_2 \widetilde{b}_2\|_{H^2}^2+\|\partial_2 \widetilde{b}_3\|_{H^2}^2\right),\\
Q_{71}\leq &C\|u\|_{H^3}\left(\|\partial_2\widetilde{u}_1\|_{H^2}^2+\|\partial_1 \widetilde{b}_1\|_{H^2}^2\right),\\
Q_{72}\leq &C\|u\|_{H^3}\left(\|\partial_1\widetilde{u}_2\|_{H^2}^2+\|\partial_1\widetilde{u}_3\|_{H^2}^2+\|\partial_2 \widetilde{b}_2\|_{H^2}^2+\|\partial_2 \widetilde{b}_3\|_{H^2}^2\right).
\end{align*}
The term $Q_4$ can be divided into four parts.
\begin{align*}
Q_4=&-\sum_{m=1}^2\sum_{i=1}^3 C_2^m\int \partial_i^m \widetilde{b}_3 \partial_3^{3-m}\overline{b}_1 \partial_i^2\widetilde{u}_1 dx-\sum_{m=1}^2\sum_{i=1}^3\sum_{k=2}^3 C_2^m\int \partial_i^m \widetilde{b}_3 \partial_3^{3-m}\overline{b}_k \partial_i^2\widetilde{u}_k dx-\sum_{i=1}^3 \int \widetilde{b}_3 \partial_3^3\overline{b}_1 \partial_i^2\widetilde{u}_1 dx\\
&-\sum_{i=1}^3\sum_{k=2}^3 \int \widetilde{b}_3 \partial_3^3\overline{b}_k \partial_i^2\widetilde{u}_k dx
:=Q_{41}+Q_{42}+Q_{43}+Q_{44}.
\end{align*}
The term $Q_{41}-Q_{44}$ can be treated by the same way of the term $N_{22}$.
\begin{align*}
Q_{41},Q_{43}\leq& C\|b\|_{H^3}\left(\|\partial_2\widetilde{u}_1\|_{H^2}^2+\|\partial_1 \widetilde{b}_3\|_{H^2}^2\right),\\
Q_{42},Q_{44}\leq &C\|b\|_{H^3}\left(\|\partial_1\widetilde{u}_2\|_{H^2}^2+\|\partial_1\widetilde{u}_3\|_{H^2}^2+\|\partial_2 \widetilde{b}_3\|_{H^2}^2\right).
\end{align*}
To bounded $Q_5$, we use the inequality (\ref{fge3}), Lemma (\ref{lempoincare}) and Young inequality to get
\begin{align*}
Q_5=&-\sum_{m=1}^2\sum_{i=1}^3 C_2^m\int \partial_i^m u\cdot\nabla \partial_i^{2-m}\widetilde{b}\cdot \partial_i^2\widetilde{b} dx\\
\leq&C\sum_{m=1}^2\sum_{i=1}^3\|\nabla \partial_i^{2-m}\widetilde{b}\|^\frac{1}{2} \|\partial_1 \nabla \partial_i^{2-m}\widetilde{b}\|^\frac{1}{2} \|\partial_i^2\widetilde{b}\|^\frac{1}{2} \|\partial_2 \partial_i^2\widetilde{b}\|^\frac{1}{2}\|\partial_i^m u\|^\frac{1}{2} \|\partial_3\partial_i^m u\|^\frac{1}{2}\\
\leq&C\|u\|_{H^2}\left(\|\partial_1 \widetilde{b}\|^2+\|\partial_2 \widetilde{b}\|^2\right).
\end{align*}
For the term $Q_6$, we just discuss the two cases.
\begin{align*}
Q_6=-\sum_{m=1}^2\sum_{i=1}^3 C_2^m\int \partial_i^m \widetilde{u}_3 \partial_3^{3-m}\overline{b}\cdot \partial_i^2\widetilde{b} dx-\sum_{i=1}^3\int \widetilde{u}_3 \partial_3^3\overline{b}\cdot \partial_i^2\widetilde{b} dx
:=Q_{61}+Q_{62}.
\end{align*}
The bound of $Q_{61}$ and $Q_{62}$ are similar to that for $N_{22}$ and $Q_{23}$, respectively.
We find
\begin{align*}
Q_{61},Q_{62}\leq& C\|b\|_{H^3}\left(\|\partial_1\widetilde{u}_3\|_{H^2}^2+\|\partial_2 \widetilde{b}\|_{H^2}^2\right).
\end{align*}
by the same way, we get
\begin{align*}
Q_8\leq& C\|u\|_{H^3}\left(\|\partial_1\widetilde{b}_3\|_{H^2}^2+\|\partial_2 \widetilde{b}\|_{H^2}^2\right).
\end{align*}
Combining these estimates of $N_1-N_8,~Q_1-Q_8$, we get
\begin{align*} 
\frac{d}{dt}&\|(\widetilde{u},\widetilde{b})\|_{H^2}^2+2\nu\|(\partial_1 \widetilde{u}_1, \partial_2 \widetilde{u}_2)\|_{H^2}^2+2\nu\left(\|(\partial_2, \partial_3) \widetilde{u}_1\|_{H^2}^2
+\|(\partial_1, \partial_3) \widetilde{u}_2\|_{H^2}^2+\|(\partial_1, \partial_2) \widetilde{u}_3\|_{H^2}^2\right)+2\mu\|(\partial_1, \partial_2) \widetilde{b}\|_{H^2}^2\\
&\leq C\|(u,b)\|_{H^3}\left(\|(\partial_2, \partial_3) \widetilde{u}_1\|_{H^2}^2
+\|(\partial_1, \partial_3) \widetilde{u}_2\|_{H^2}^2+\|(\partial_1, \partial_2) \widetilde{u}_3\|_{H^2}^2+\|(\partial_1, \partial_2) \widetilde{b}\|_{H^2}^2\right).
\end{align*}
We now specify $\lambda=min\{\nu, \mu\}$. Let $\epsilon>0$ sufficiently small and $\|(u_{in},b_{in})\|_{H^3} \leq\epsilon$ under the assumptions of Theorem (\ref{theorem1t}). Then
$$\|(u,b)\|_{H^3}\leq\epsilon,~ 2\lambda-\|(u,b)\|_{H^3}\geq \lambda.$$
By Lemma (\ref{lempoincare}), we finally obtain that
\begin{align*} 
\frac{d}{dt}&\|(\widetilde{u},\widetilde{b})\|_{H^2}^2+\lambda\|(\widetilde{u},\widetilde{b})\|_{H^2}^2\leq 0.
\end{align*}
From Gr\"{o}nwall inequality, we get the exponential decay result
\begin{align*}
\|(\widetilde{u}, \widetilde{b})\|_{H^2}\leq \|(u_{in}, b_{in})\|_{H^2}e^{-\frac{\lambda}{2} t}.
\end{align*}
This ends the proof of Theorem \ref{theorem2t}.
\end{proof}

\noindent\textbf{Acknowledgments.} The work of Aibin Zang was partially supported by the National Natural Science Foundation of China (Grant no. 12261039, 12061080), and Jiangxi Provincial Natural Science Foundation  (No. 20224ACB201004). Yuelong Xiao was supported in part by National Natural Science Foundation of China (Grant no.11871412). The research of Xuemin Deng partially Supported by Postgraduate Scientific Research Innovation Project of Hunan Province (Grant no.CX20210603, XDCX2021B096).


\end{document}